\documentclass[12pt,twoside]{amsart}
\usepackage{amsmath}
\usepackage{amsfonts}
\usepackage{amssymb}
\usepackage{amsthm}
\usepackage[english]{babel}

\newtheorem{theo}{Theorem}[section]
\newtheorem{lem}[theo]{Lemma}
\newtheorem{prop}[theo]{Proposition}
\newtheorem{cor}[theo]{Corollary}
\newtheorem{defn}{Definition}[section]
\newtheorem{claim}{Claim}

\newcommand{\N}{{\mathbb N}}

\renewcommand{\S}{{\mathcal S}}
\newcommand{\C}{{\mathbb C}}
\renewcommand{\SS}{{\mathbb S}}
\renewcommand{\H}{${\bf (H) } $} 
\newcommand{\Ht}{${\bf (\tilde H) } $}
\newcommand{\RR}{{\mathcal R}}
\newcommand{\HH}{{\mathcal H}}

\def\Z{\mathbb{Z}}

\def\R{\mathbb{R}}

\begin{document}

\title
[Convolution powers]
{ On the Ritt property and weak type maximal inequalities for 
convolution powers on $\ell^1(\Z)$}

%\end{center}

\author{Christophe Cuny }
\address{Laboratoire MICS, Centralesupelec, Grande Voie des Vignes, 92295 Chatenay-Malabry cedex, FRANCE.}
%\newline\indent {\it current address:}
%Laboratoire MAS, Ecole Centrale de Paris}
\email{christophe.cuny@ecp.fr}
\thanks{I am very grateful to Alexander Gomilko and, more 
particularly, to Yuri Tomilov 
who  both noticed several inaccuracies in a previous version. 
The paper substantially benefited from our discussions.}

%\subjclass{Primary: 47A35, 28D05, 37A05; Secondary: 47B38}

\begin{abstract}
In this paper we study the behaviour of convolution powers of probability 
measures $\mu$ on $\Z$, such that $(\mu(n))_{n\in \N}$ is completely monotone or such that $\nu$ is centered  with a second moment.
In particular we exhibit many new examples of probability measures on $\Z$ 
having the so called Ritt property and whose convolution powers satisfy 
weak type maximal inequalities in $\ell^1(\Z)$.   
\end{abstract}

\maketitle

%We prove a new maximal inequality for (adapted) stationary processes taking 
%values in a separable real Banach space. In the case of smooth Banach spaces, we derive the Marcinkiewicz-Zygmund strong law 
%of large numbers, the weak invariance principle and the almost sure 
%invariance principle under conditions similar to a condition introduced 
%by Maxwell and Woodroofe. Except for the weak invariance principle our 
%results are also new  in the real case. The results extend to non-adapted 
%processes and partly to processes arising in non-invertible dynamical systemes.


\section{Introduction}

Let $\mu$ be a probability on $\Z$. Given an invertible 
bi-measurable transformation $\tau$ on a measure space 
$(\SS,\S,\lambda)$ we define a \emph{positive contraction} of every $L^p(\lambda)$, 
$1\le p\le \infty$, by setting 
$$
P_\mu(\tau)(f):= \sum_{k\in \Z} \mu(k)f\circ \tau^k \qquad \forall f\in L^p(m)
\, .
$$

Several authors, see for instance \cite{BJR90}, \cite{BJR92}, \cite{Reinhold},
\cite{BJR94}, \cite{JR}, \cite{BC}, \cite{Losert01}, \cite{Wedry}, 
\cite{RSW}, \cite{RS} or \cite{CCLritt}, studied the almost everywhere behaviour of the iterates 
of $\mu(\tau)$, i.e. of $(\mu^{*n}(\tau))_{n\ge 1}$, acting on $L^p(\lambda)$, $1\le p<\infty$. 

\medskip

When $p>1$, the almost everywhere behaviour has been characterized thanks to the so called \emph{bounded angular ratio} property,  introduced in 
\cite{BJR94} and which is  equivalent to the so called \emph{Ritt} property on 
$\ell^p(\Z)$, $p>1$. Let us recall the definition of those properties.

\smallskip

\begin{defn}
Let $\mu$ be a probability measure on $\Z$. We say that $\mu$ is 
strictly aperiodic if $|\hat \mu(\theta)|<1$ for every $\theta\in (0,2\pi)$. 
We say that $\mu$ has bounded angular ratio (BAR) if moreover 
\begin{equation}\label{BAR}
\sup_{\theta\in (0,2\pi)}\frac{|1-\hat \mu(\theta)|}{1-|\hat \mu(\theta)|}
<\infty\, .
\end{equation}
\end{defn}

The strict aperiodicity is equivalent to the fact that the support of 
$\mu$ is not contained in a coset of a proper subgroup of $\Z$. In particular, it holds 
whenever the support of $\mu$ contains two consecutive integers.

\begin{defn}\label{defritt}
We say that a probability measure $\mu$ on $\Z$ is Ritt on $\ell^p(\Z)$, for some 
$p\ge 1$, if 
$$\sup_{n\ge 1}n\|\mu^{*n}-\mu^{*(n+1)}\|_{\ell^p(\Z)} 
<\infty\, .
$$
When $p=1$ we say simply that $\mu$ is Ritt, because then, it is Ritt 
on all $\ell^r(\Z)$, $r\ge 1$.\\
Denote by $\RR$ the set of Ritt probability measures on $\Z$.
\end{defn}

\medskip

A version of the next Theorem may be found for instance in Cohen, Cuny and Lin \cite[Theorem 4.3]{CCLritt}. Their Theorem 4.3 is not formulated exactly as below but the proof of Theorem \ref{CCL} may be done  similarly.  
The equivalence the item $(vi)$ with the other items follow from 
their Proposition 6.4. In all the paper we use the notation $\N:=\{0,1,2\ldots\}$.

\begin{theo}\label{CCL}
Let $\mu$ be a strictly aperiodic probability on $\Z$. The following are equivalent: 
\begin{itemize}
\item [$i)$] $\mu$ has BAR;
\item [$ii)$] There exist $p>1$ and $C_p>0$ such that for every invertible 
bi-measurable transformation $\tau$ on a measure space $(S,\S,\lambda)$, 
\begin{equation}\label{inemaxritt}
\|\sup_{n\ge 1}|(P_\mu(\tau))^nf|\,  \|_{p,\lambda}\le C_p \|f\|_{p,\lambda}\qquad \forall f\in L^p(\lambda)\, ;
\end{equation}
\item [$(iii)$] There exists  $p>1$ such that for every invertible 
bi-measurable transformation $\tau$ on a probability space $(S,\S,\lambda)$ and every 
$f\in L^p(\lambda)$, $((P_\mu(\tau))^nf)_{n\in \N}$ converges $\lambda$-a.e.
\item [$(iv)$] There exist $p>1$ and $C_p>0$ such that  
\begin{equation}\label{inemaxshift}
\|\sup_{n\ge 1}|(P_\mu(R))^{n}f|\,  \|_{\ell^p(\Z)}\le C_p \|f\|_{\ell^p(\Z)}
\qquad f\in \ell^p(\Z)\, ,
\end{equation}
where $R$ is the right shift on $\Z$;
\item [$v)$] There exists $p>1$ such that $\mu$ is Ritt on ${\ell^p(\Z)}$.
\item [$(vi)$] There exists $p>1$ such that for every $m\in \N$, there 
exists $C_{m,p}>0$ such that 
\begin{equation*}
\|\sup_{n\ge 1}n^m |(I-P_\mu)^m(P_\mu(R))^{n}f|\,  \|_{\ell^p(\Z)}\le 
C_{m,p} \|f\|_{\ell^p(\Z)}
\qquad f\in \ell^p(\Z)\, ,
\end{equation*}
\end{itemize}
\end{theo}

Actually,  if any of the above properties  
holds, then the conclusion of  $(ii)$, $(iii)$, $(iv)$ or 
$(v)$ holds for all $p>1$. 

\smallskip

The proof of the above theorem follows from recent works of Le Merdy and Xu,  \cite{LX12} and \cite{LX13}, who studied positive Ritt contractions $T$ of $L^p(\SS,\S,m)$ ($p$ being fixed). Recall that a contraction $T$ on a Banach space $X$ is Ritt if  
 $\sup_{n\in \N}n\|T^n-T^{n+1}\|_X<\infty$, which is compatible with our definition \ref{defritt} which just says that the operator of convolution 
 by $\mu$ is Ritt  on $X=\ell^p(\Z)$.

 Le Merdy and Xu proved that any positive Ritt contraction satisfies maximal inequalities in spirit of \eqref{inemaxritt}. They also obtained square function estimates, oscillation inequalities and variation inequalities. See also 
\cite{CCLritt} for related results. 

%Of course, once one has $ii)$, one has the $m$-a.e. convergence 

In this paper we are concerned with the case when $p=1$, and we address the following two questions. 

\smallskip

{\bf - Question 1:} For what probability measures $\mu$ on $\Z$ does one have a weak type $(1,1)$-maximal inequality: 
\begin{equation}\label{weakL1}
 \#\{k\in \Z \, :\, \sup_{n\ge 1}|\mu^{*n}*f|>\lambda\} 
\le \frac{C}\lambda \|f\|_{\ell^1(\Z)} \qquad \forall \lambda\ge 0\, .
\end{equation}
More generally, given $m\in \N$, does there exist $C_m>0$ such that 
(with the convention $(\delta_0-\mu)^{*0}=\delta_0$)
 \begin{equation}\label{bc}
 \sup_{\lambda >0}\lambda  \#\{k\in \Z\, :\, \sup_{n\ge 1} 
n^m |\mu^{*n}*(\delta_0-\mu)^{*m}*f(k)|\ge \lambda \} \le C_m 
\|f\|_{\ell^1(\Z)}\qquad \forall f\in \ell^1(\Z)\, .
 \end{equation}

{\bf - Question 2:} For what probability measures $\mu$ on $\Z$ does one have 
the Ritt property in $\ell^1(\Z)$: 
\begin{equation}\label{Ritt}
\sup_{n\ge 1}n\|\mu^{*n}-\mu^{*(n+1)}\|_{\ell^1(\Z)} <\infty\, .
\end{equation}

%In this paper, we wish to produce new examples of probabilities on $\Z$ with 
%the BAR property. Moreover, for those examples we prove a weak type 
%$(1,1)$ version of \eqref{inemaxshift}

Notice that if $\mu$ satisfies \eqref{weakL1} then, by the Marcinkiewicz interpolation theorem (between weak $L^1$ and $L^\infty$), it does satisfy \eqref{inemaxritt}, hence $\mu$ has BAR. Notice also that if $\mu$ satisfies \eqref{Ritt} then, by Theorem 
\ref{CCL} has also BAR. Hence, the questions we intend to answer are:  what extra conditions, in addition to the BAR property, 
are sufficient to have 
\eqref{weakL1}, \eqref{bc} or \eqref{Ritt} ? 

\smallskip

Let us discuss the known results concerning those questions, before presenting our results. As far as we know, when $m\ge 1$, \eqref{bc} has not been investigated 
before.

\smallskip

The simplest examples of probability measures having BAR are the symmetric 
ones. Bellow, Jones and Rosenblatt \cite{BJR94} proved that if $\mu$ 
is symmetric such that $(\mu(n))_{n\ge 0}$ is non increasing then 
\eqref{weakL1} holds. We do not know whether  \eqref{Ritt} 
holds as well, in this case, but we provide sufficient conditions 
in Section \ref{symmetric}.

Another case where \eqref{weakL1} holds is when $\sum_{k\in \Z} 
k^2 \mu(k) <\infty$ ($\mu$ has a second moment) and 
$\sum_{k\in \Z} k\mu(k)=0$. This has been proved by Bellow and Calder\'on 
\cite{BC}. Again the Ritt property is not known in that case. The proof 
of Bellow and Calder\'on is based on general intermediary results that have been 
extended recently by Wedrychowicz \cite{Wedry}. Wedrychowicz proved that 
\eqref{weakL1} holds for centered probability measures (hence with a first moment)   
having BAR and satisfying some extra conditions. Examples without second moment 
are presented in \cite{Wedry}.

\medskip

Several examples of probabilities having the Ritt property in $\ell^1(\Z)$ 
may be found in Dungey \cite{Dungey11}, see sections 4 and 5 there.

\medskip

Let us now present our results. As mentionned above, the method 
of Bellow and Calder\'on is fairly general. Actually, if one follows carefully their paper, one realizes that the following definition comes somewhat naturally into play. 

\begin{defn}
We say that a probability measure $\mu$ on $\Z$ satisfies the hypothesis ${\bf (H)}$ if 
$\hat \mu$ is twice continuously differentiable on $[-\pi,\pi]-\{0\}$ and if there exists an even and continuous function $\psi$ on $[-\pi,\pi]$, vanishing at $0$ 
and continuously differentiable  on $[-\pi,\pi]-\{0\}$, and 
some constants $c,C>0$ such that for every $\theta\in (0,\pi]$
\begin{itemize}
\item [$(i)$] $|\hat \mu (\theta)|\le 1-c\psi(\theta)$;  
\item [$(ii)$]  $|\theta\hat \mu '(\theta)|\le C \psi(\theta)$;
\item [$(iii)$] $| \hat\mu '(\theta)|\le C \psi'(\theta)$; 
\item [$(iv)$] $|\theta\hat \mu''(\theta)|\le C \psi'(\theta)$.
\end{itemize} 
Let us denote by $\HH$ the set of probability measures satisfying hypothesis 
${\bf (H)}$.
\end{defn}

The relevance of the hypothesis ${\bf (H)}$ lies in the following, where 
we also give stability properties of $\HH$ as well as of $\RR$. We say that 
a set of probability measures on $\Z$ is stable by symmetrization if 
whenever $\mu=(\mu(n))_{n\in \Z}$ belongs to that set so does 
$\check \mu=(\mu(-n))_{n\in \Z}$.

\begin{theo}\label{theo-H}
~~\begin{itemize}
\item [$(i)$] The set $\HH$ is convex and  stable by convolution and by symmetrization. 
\item [$(ii)$] The set $\RR$ is convex and stable by convolution and by symmetrization.  
\item [$(iii)$] Let $\mu\in \HH$. Then, $\mu$ satisfies \eqref{weakL1}.
\item [$(iv)$] Let $\mu\in \HH \cap \RR$. Then, for every 
$m\in \N$, there exists $C_m>0$ such that $\mu$ satisfies 
\eqref{bc}.
\end{itemize}
\end{theo}

Theorem \ref{theo-H} follows from several results: item $(i)$ 
follows from Proposition \ref{ext}, item $(ii)$ may be proved as Proposition 
\ref{prop-conv}, items $(iii)$ and $(iv)$ follow from 
Proposition \ref{genprop}.

\medskip

Then, our goal is to provide many examples of elements of $\HH\cap \RR$.

%In this paper, we intend to exhibit large classes of examples of Ritt %probability measures satisfying ${\bf (H)}$, hence for which 
%Theorem \ref{theo-H} applies. Moreover, we shall see that the set of probability measures satisfying  hypothesis ${\bf (H)}$ is convex and 
%stable by convolution. On the other hand, it is known that the set of Ritt probability 
%measures is also convex and stable by convolution. 

Our first examples are the ones already considered by Bellow and Calder\'on, 
in particular the fact that $\mu$ as in the next theorem satisfies \eqref{weakL1} is not new, while the Ritt property is new.  The proof of Theorem \ref{centered} is done in Section 2.3.

\begin{theo}\label{centered}
Let $\mu$ be a centered and strictly aperiodic probability 
measure on $\Z$ with finite 
second moment. Then $\mu\in \HH \cap \RR$.
\end{theo}

Then, we shall consider probability measures 
$\mu$, such that $(\mu(n))_{n\ge 0}$ is completely monotone 
(see the next section for the definition). In this context we are able to 
 characterize  the BAR property.  The idea of considering completely monotone sequences 
 was motivated by Gomilko-Haase-Tomilov \cite{GHT} and 
 Cohen-Cuny-Lin \cite{CCLritt}.
 
 \smallskip
 
 \begin{theo}\label{theoNC}
 Let $\mu$ be a probability measure on $\Z$ supported on $\N$, such that 
 $(\mu(n))_{n\in \N}$ is completely monotone. Then 
 \begin{itemize} 
 \item [$(i)$] $\mu$ has BAR 
 if and only if there exists $C>0$ such that 
 \begin{equation}\label{cns1}
 \sum_{k=1}^n k\mu(k)\le C n\sum_{k\ge n}\mu(k)\qquad \forall n\ge 1\, .
 \end{equation}
 \item [$(ii)$] Assume that $\mu$ has BAR.  Let $\sigma$ be a probability measure 
  on $\Z$ such that $\sum_{n\in \Z}n^2\sigma(n)<\infty$. Then, 
  $\mu *\sigma\in \HH\cap \RR$ and for every 
 $\alpha\in (0,1]$ $\alpha\mu +(1-\alpha)\sigma\in \HH\cap \RR$. 
 In particular (take $\sigma=\delta_0$), $\mu \in \HH\cap \RR$.
 \end{itemize}
 \end{theo}
 \noindent {\bf Remarks.}  Notice  that we do not assume 
 $\sigma$ to be centered. The conclusion of item $(ii)$ actually 
 holds for $\sigma$ such that $\hat \sigma$ is twice continuously 
 differentiable on 
 $[-\pi,\pi]-\{0\}$ with $\hat \sigma'$ and $\theta\mapsto \theta 
 \hat \sigma''(\theta)$  bounded. Moreover, see Proposition \ref{prop-conv}, 
 it is possible to relax the conditions on $\hat \sigma$ if one is only concerned with the Ritt Property. We were not able to provide a perturbation result in the spirit 
 of Theorem \ref{theocentered}.
 \medskip
 
Item $(i)$ follows from Propositions \ref{prop} and \ref{propcarac}. 
Item $(ii)$ is proved in sections 3.2 and 3.3. 

\smallskip 

The proof of the Ritt property in Theorem \ref{theoNC} is based on a recent of Gomilko and Tomilov \cite{GT}. The fact that when $\mu$ has BAR 
 $\delta_1*\mu$ is Ritt has been proven by Gomilko and Tomilov 
 \cite{GT1}, see their Theorem 7.1. Their proof is also based on \cite{GT}.
 
 %Our main results are the following (see section 3 for the definition of 
 %completely monotone sequences).  In the next theorem, the fact that when $\mu$ has 
 %BAR it is  Ritt has been proved by Gomilko and Tomilov \cite{GT1}. 
 %The main argument makes use of their previous work \cite{GT}, hence 
 %the result should be attributed to Gomilko and Tomilov. However, \eqref{bc} and the characterization of the BAR property by \eqref{cns1}
 %are  new.

 %Actually, \eqref{cns1}, \eqref{ritt} and \eqref{bc} 
 %are all equivalent to the fact that $\mu$ has BAR. Notice that \eqref{cns1} implies that 
 %$\mu$ does not have a first moment (i.e. $\sum_{k\ge 1}k\mu(k)=+\infty$). 
 %The proof of \eqref{ritt} makes use of a difficult result of Gomilko and Tomilov 
 %\cite{GT}. 
 %If \eqref{ritt} holds then (see the proof of Proposition \ref{genprop}) for every $m\ge 1$, 
 %\begin{equation}
%\sup_{n\in \N} n^m\|\mu^{*n}(\delta_0-\mu)^{*m}\|_{\ell^1(\Z)}<\infty\, , 
% \end{equation}
 %The condition \eqref{cns1} (as well as the condition \eqref{cns2}) 
 %is related to the work of Losert \cite{Losert99}, see Propositions 1,3 and 4 there.
  \begin{theo}\label{theocentered}
 Let $\mu$ be a \emph{centered} probability measure on 
 $\Z$ supported on $\{-1\}\cup\Z$, such that 
 $(\mu(n))_{n\in \N}$ is completely monotone. Then 
 \begin{itemize}
 \item [$(i)$] $\mu$ has BAR 
 if and only if there exists $C>0$ such that 
 \begin{equation}\label{cns2}
 n\sum_{k\ge n} k\mu(k)\le C \sum_{k=1}^n k^2\mu(k)\qquad 
 \forall n\ge 1\, .
 \end{equation}
\item [$(ii)$] Assume that $\mu$ has BAR. Let $\sigma$ be a centered 
probability measure on 
 $\Z$, such that there exists $c>0$ such that 
 $\sum_{n\in \Z}n^2|\sigma(n)-c\mu(n)|<\infty$. Then, $\sigma\in \HH\cap \RR$. In particular (take $c=1$ and $\sigma=\mu$), $\mu\in \HH\cap \RR$.
 
 \end{itemize}
 \end{theo}
% \noindent {\bf Remark.} Taking $c=1$ and $\sigma=\mu$ in item $(ii)$ we see that $\mu\in \HH\cap \RR$. 

\smallskip

Moreover (see section 5), we also study symmetric probability measures with completely monotone coefficients.

%In Theorem \ref{centered} we do not assume that $(\mu(n))_{n\in \N}$ 
%be completely monotone.
% Moreover, if $\mu$ has BAR, \eqref{ritt} holds and for every $m\in \N$, 
% there exists $C_m>0$ such that \eqref{bc} holds.
% \end{theo}
% \noindent {\bf Remark.} Notice that if $\mu$ has a moment of order 2, 
% then \eqref{cns2} is satisfied. 

% Moreover, we prove that if such a 
% $\mu$ has BAR then \eqref{Ritt} and \eqref{weakL1} hold as well. Then, we prove 
% othe weak type maximal inequalities providing a speed of convergence 
% in \eqref{weakL1} when $f=(\delta_0-\mu)^{*k}*g$ for some $g\in \ell^1(\Z)$. 

In the above theorems, we obtain weak type maximal inequalities in 
$\ell^1(\Z)$. Of course, by mean of  transference principles (see e.g. 
\cite{Bellow} or \cite[page 164]{Weber}, one may derive similar results for the operator $P_\mu(\tau)$ in the spirit of 
\eqref{inemaxritt} as well as some almost everywhere convergence results. 
We leave that "standard" task to the reader.

\medskip 

The paper is organized as follows. In Section 2, we prove Theorem 
\ref{theo-H} and prove the Ritt property under a slightly stronger 
assumption than hypothesis ${\bf (H)}$.  In section 3, we consider probability measures as in  Theorem \ref{theoNC} and prove Theorem \ref{theoNC}
In section 4, we consider probability measures as in Theorem \ref{theocentered} and prove Theorem \ref{theocentered}. 
In section 5 we consider symmetric probability measures. Finally, in section 6 
we discuss several open questions on the topic.

\medskip

Before going to the proofs, we would like to mention that the above theorems provide new situations to which the results of Cuny and Lin \cite{CL} apply, see 
examples 1 and 2 there. 

\section{General criteria for  maximal inequalities and for the Ritt property}

In this section we give general conditions ensuring 
weak type maximal inequalities associated with sequences of probabilities on $\Z$ 
and conditions ensuring the Ritt property. 

In the case of weak type maximal inequalities, the obtained conditions are derived from slight modifications of   known results, see e.g. \cite{BC} and 
\cite{Wedry}.

\medskip

\subsection{Sufficient conditions for weak type maximal inequalities}

We start with  the following result of Bellow and 
Calder\'on \cite{BC}, see also Zo \cite{zo} for a related result. 
Actually, Bellow and Calder\'on considered only the case of probability measures, but their proof extends to the situation below.

\begin{theo}[Bellow-Calder\'on]\label{theoBC}
Let $(\sigma_n)_{n\in \N}$ be a sequence of finite \emph{signed} measures  on $\Z$ such that $\sup_{n\in \N} \|\sigma_n\|_{\ell^1}<\infty$. 
Assume that there exists $C>0$ such that $k,\ell\in \Z$ 
 with $0< 2|k|\le \ell$, 
 \begin{equation}\label{BCcond}
 |\sigma_n(k+\ell)-\sigma_n(\ell)|\le C\frac{k}{\ell^2}\qquad 
 \forall n\in \N\, .
 \end{equation}
 Then, there exists $C'>0$ such that for every $f\in \ell^1(\Z)$, 
\begin{equation*}
\sup_{\lambda >0} \#\{k\in \Z\,:\, \sup_{n\in \N} |\sigma_n*f(k)|\ge \lambda\}
\le \frac{C'}{\lambda}\|f\|_{\ell^1}\, .
\end{equation*}
\end{theo}

In order to apply Theorem \ref{theoBC} we shall need the following version 
of  Corollary 3.4 of \cite{BC}. 

\begin{lem} \label{BClemma}
Let $(\sigma_n)_{n\in \N}$ be a sequence of finite signed 
measures, such that for every $n\in \N$, $\hat \sigma_n$ is twice continuously differentiable on $\R-2\pi\Z$. If moreover 
\begin{equation}\label{BCcond2}
\sup_{n\in \N} \int_{-\pi}^\pi |\theta \hat \sigma_n''(\theta)|
d\theta <\infty \, ,
\end{equation}
and 
\begin{equation}\label{limit}
\lim_{\theta\to 0,\theta\neq 0} \theta\hat \sigma_n'(\theta)=0\qquad \forall n\in \N\, ,
\end{equation}
then \eqref{BCcond} holds. 
\end{lem}
\noindent {\bf Remark.} It  follows from \eqref{BCcond2} 
(and the continuity of $\hat \sigma_n$ at 0) that the limit in \eqref{limit} 
exists, hence condition \eqref{limit} is just that the limit is 0. 

\noindent {\bf Proof.} For every $k\in \Z-\{0\}$, we have
$\sigma_n(k)=\int_{-\pi}^\pi \hat \sigma_n(\theta)
{\rm e}^{-ik\theta}d\theta$. 
Let $\pi >\varepsilon>0$. Performing two integration by parts 
as in \cite{BC} to evaluate $\int_{\varepsilon}^\pi \hat \sigma_n(\theta){\rm e}
^{-ik\theta}d\theta$ and $\int_{-\pi}^{-\varepsilon} \hat \sigma_n(\theta){\rm e}
^{-ik\theta}d\theta$, using our assumptions and letting $\varepsilon
\to 0$, we see that 
$$
\sigma_n(k)=\int_{-\pi}^\pi \hat \sigma_n''(\theta)\frac{1-{\rm e}
^{-ik\theta}}{k^2}\, d\theta\, .
$$
Then, we conclude as in \cite{BC}. \hfill $\square$

%We first give a general result ensuring that condition 
%\eqref{BCcond} be satisfied. 

\medskip

\begin{prop}\label{genprop}
Let $\mu$ be a probability measure on $\Z$  satisfying 
 hypothesis ${(\bf H)}$. Then,  for every $m\in \N$, 
\begin{equation}\label{BCcond3}
\sup_{n\ge 1}n^m\int_{-\pi}^\pi |\theta||(\hat \mu^n(1-\hat \mu)^m)''(\theta)|
 d\theta <\infty\, ,
 \end{equation} 
 and 
 \begin{equation*}
\lim_{\theta\to 0,\theta\neq 0} \theta(\hat \mu^n(1-\hat \mu)^m)'(\theta)=0\qquad \forall n\in \N\, ,
\end{equation*}
 In particular, there exists $C>0$ such that for every $f\in \ell^1(\Z)$, 
\begin{equation}\label{weakinemax0}
\sup_{\lambda >0} \lambda \#\{k\in \Z\,:\, \sup_{n\ge 1}|\mu^{*n}*f(k)|\ge \lambda\}\le C\|f\|_{\ell^1}\, .
\end{equation}
 If moreover $\mu$ is Ritt then, for every $m\ge 1$,
there exists $C_m>0$ such that for every $f\in \ell^1(\Z)$, 
\begin{equation}\label{weakinemax}
\sup_{\lambda >0} \lambda \#\{k\in \Z\,:\, \sup_{n\ge 1} n^m|\mu^{*n}*(\delta_0-\mu)^{*m}*f(k)|\ge \lambda\}
\le C_m\|f\|_{\ell^1}\, .
\end{equation}
\end{prop}
\noindent {\bf Remarks.} The proposition is related to Theorem 2.10 of 
Wedrychowicz \cite{Wedry}. Notice that, by $(ii)$, $\psi$ is non-negative and by $(iii)$ it is non-decreasing. We shall see in proposition \ref{genprop2} 
that if there exists $C>0$ such that for every $\theta\in (0,\pi]$, 
$\psi(\theta)\le C\theta \psi'(\theta)$, then $\mu$ is automatically Ritt. 
\\
\noindent {\bf Proof. } If $\mu=\delta_0$ the result is trivial. 
Hence we assume that $\mu\neq \delta_0$. In particular, by $(ii)$, 
$\psi$ cannot vanish in a neighbourghood of $0$, hence is positive on 
$(0,\pi]$. Then $|\hat \mu|<1$ on $(0,\pi]$ (hence $\mu$ is strictly aperiodic).

\smallskip
 
 We have, on $(0,\pi]$. 
\begin{gather}\label{mu''}
n^m|(\hat \mu^n(1-\hat \mu)^m)''|\le n^{m+2}|\hat \mu|^{n-2}|\hat \mu'|^2
|1-\hat \mu|^m + 
2m n^{m+1} |\hat \mu|^{n-1}|\hat \mu'|^2|1-\hat \mu|^{m-1}\\
\nonumber \qquad  +n^{m+1}|\hat \mu|^{n-1} |\hat \mu''|
|1-\hat \mu|^{m-1} + mn^m |\hat \mu|^n |\hat \mu''|\, |1-\hat \mu|^{m-1} 
+m(m-1)n^m|\hat \mu |^n |\hat \mu'|^2|1-\hat \mu|^{m-2} \, . 
\end{gather}
Using $(i)$ and the fact that $\psi$ is continuous with $\psi(0)=0$, there exist $\eta>0$ and $c >0$, such 
that for every $\theta\in [0,\eta]$, 
\begin{equation}\label{estmu}
|\hat \mu(\theta)|\le {\rm e}^{-c\psi(\theta)}\, .
\end{equation}
Since, $\sup_{\theta\in [\eta,\pi]} |\hat \mu(\theta)|<1$, 
taking $c$ smaller if necessary, we may assume that \eqref{estmu} holds 
for every $\theta\in [0,\pi]$. 

\smallskip

Using $(iii)$ and that $\psi$ is continuous at 0, we see that 
$\psi'$ and $\hat\mu'$ are in $L^1$, hence that 
\begin{equation}\label{estimate}
|1-\hat \mu(\theta)|\le \psi(\theta)\qquad \forall \theta\in (0,\pi]
\, . 
\end{equation}
\smallskip

Combining \eqref{estmu} and \eqref{estimate} with $(ii)$, $(iii)$ and $(iv)$ (taking care with the cases $m=0$ and $m=1$),  
we see that there exists $C_m>0$ such that for every 
$\theta\in (0,\pi]$ and every $n\in \N$, 
\begin{equation*}
|\theta (\hat\mu^n)''(\theta)|\le C_m\big( n^2\psi^{m+1}(\theta)+
n\psi^{m}(\theta)+ (m-1)\psi^{m-1}(\theta)\big){\rm e}^{-c(n-2)\psi(\theta)}\, n^m\psi'(\theta)\, .
\end{equation*}

Using that the integrand below is even and the change of variable 
$u=(n-2)\psi(\theta)$, we see that 
\begin{gather*}
\sup_{n\in \N} \int_{-\pi}^\pi |\theta||(\hat \mu^n(1-\hat \mu)^m)''(\theta)|
 d\theta \\ \le \tilde C_m \int_0^{+\infty} (u^{m+1}+u^m 
 +(m-1)u^{m-1}){\rm e}^{-cu} \, du \, <\infty,
\end{gather*}
and \eqref{BCcond3} holds. 

\smallskip

The fact that \eqref{limit} holds follows from item $(ii)$, using 
 that $\psi$ is continuous at 0, with $\psi(0)=0$, and that $\hat \mu$ is 
 bounded.

\smallskip

 Let $m\in \N$ and set $\sigma_n=\sigma_{n,m}:= 
n^m \mu^{*n}(\delta_0-\mu)^{*m}$ for every $n\in \N$.  

\smallskip

It follows from Theorem 
\ref{theoBC} and Lemma \ref{BClemma}, that \eqref{weakinemax0} holds. 

When $m\ge 1$, it follows from Theorem 
\ref{theoBC} and Lemma \ref{BClemma}, that \eqref{weakinemax} holds 
provided that
\begin{equation}\label{rittcond}
\sup_{n\in \N} \|\sigma_{n,m}\|_{\ell^1(\Z)}<\infty\, .
\end{equation}

When $m=1$, \eqref{rittcond} is just the definition of the Ritt property. 

Let $m\ge 2$, write $n=m\ell +k $, with $\ell\in \N$ and $0\le k\le m-1$. 
We have 
$$
\|n^m \mu^{*n}*(\delta_0-\mu)^{*m}\|_{\ell^1(\Z)} 
\le  m^m \|(\ell +1) \mu^{*\ell}*(\delta_0-\mu)\|_{\ell^1(\Z)}^m \,, 
$$
and the latter is bounded uniformly with respect to $\ell\in \N$, 
since $\mu$ is Ritt. \hfill $\square$

\smallskip

To conclude this subsection we shall study stability properties of set of probabilities satisfying the weak type maximal inequalities. 

\smallskip 

It is well-known, see e.g. Proposition 3.2 of \cite{Dungey11}, that the set 
of Ritt probability measures on $\Z$  is convex and stable by convolution. 
Actually, \cite{Dungey11} deals with probability measures supported by $\N$, 
but the proof is the same. 

Let $p>1$. It is not difficult to see that the set of probability measures $\mu$ 
on $\Z$, such that there exists $C_p>0$ such that \eqref{inemaxshift} 
holds is also convex and stable by convolution. 

However it is unclear (and probably not true) whether the 
set of probability measures $\mu$ on $\Z$ satisfying \eqref{weakinemax} 
for every $m\in \N$ (or for some $m\in \N$) is also convex and stable by convolution. 
Nevertheless, we have the following. 

\begin{prop}\label{ext}
Let $\mu_1$ and $\mu_2$ be probability measures satisfying 
hypothesis ${\bf (H)}$. Let $\alpha \in (0,1)$. Then, $\check \mu_1$, $\mu_1
*\mu_2$ and $\alpha \mu_1+(1-\alpha)\mu_2$ satisfy hypothesis ${\bf (H)}$.
\end{prop}
\noindent {\bf Remark.} Recall that $\check \mu_1$ is the probability measure defined by $\check \mu_1(n)=\mu_1(-n)$ for every $n\in \Z$. \\
\noindent {\bf Proof.} The fact that $\check \mu_1$ satisfies
hypothesis ${\bf (H)}$ is obvious.

Let $\psi_i$, $c_i,C_i$ be the terms associated 
with $\mu_i$ ($i\in \{0,1\}$) such that the items $(i)-(iv)$ 
of hypothesis \H\,  be satisfied. 

Define $\mu:=\mu_1*\mu_2$ and $\psi:=  c_1 \psi_1+
c_2\psi_2$. Let $\theta\in (0,\pi]$. We have 
\begin{gather*}
|\hat \mu (\theta)|= |\hat \mu_1 (\theta)|\, |\hat \mu_2 (\theta)| 
\le 1- \psi(\theta)+c_1c_2\psi_1(\theta)\psi_2(\theta). 
\end{gather*}
Since $\psi_1$ and $\psi_2$ are continuous with $\psi_1(0)=\psi_2(0)=0$, 
there exist $c\in (0,1)$ and $\eta\in (0,\pi)$ such that 
$c_1c_2\psi_1(\theta)\psi_2(\theta)\le (1-c)\psi(\theta)$ for every 
$\theta \in (0,\eta)$. 

Hence, $|\hat \mu|\le 1-c\psi$ on $(0,\eta)$. Arguing as in the previous proof, we see that taking $c$ smaller if necessary, the inequality holds 
on $(0,\pi]$ either. 

\smallskip

Using that $\hat\mu'=\hat\mu_1'\hat \mu_2+\hat\mu_1\hat \mu_2'$, we infer 
that items $(ii)$ and $(iii)$ of hypothesis \H\,  hold. 
\smallskip

We have $\hat\mu''=\hat\mu_1''\hat \mu_2+2\hat\mu_1'\hat \mu_2'+\hat\mu_1\hat \mu_2''$. Hence, for every $\theta\in (0,\pi]$
$$
|\theta \hat\mu''(\theta)|\le C_1\psi_1'(\theta) + 2C_1\psi_1(\theta)
C_2\psi'_2(\theta) +C_2\psi_2'(\theta)\, ,
$$
and we see that item $(iv)$ holds, since $\psi_1$ is bounded. 

\medskip

Let $\alpha\in (0,1)$. Let $\mu:=\alpha\mu_1+(1-\alpha)\mu_2$. 
One can see that items $(i)-(iv)$ of hypothesis \H \, hold 
with $\psi:= \alpha  \psi_1+(1-\alpha)\psi_2$. \hfill $\square$

\subsection{A sufficient condition for the Ritt property}

In this subsection we derive a condition ensuring that a probability 
measure is Ritt. This condition will be used for centered probability measure 
with either a second moment, or a first moment and completly monotone 
coefficients. For non centered probability measure another argument will be needed. 

\smallskip

We start with a general result.

\begin{prop}\label{genprop2}
Let $(\sigma_n)_{n\in \N }$ be a sequence of finite signed measures on $\Z$, such that for every $n\in \N$, $\hat \sigma_n$ is twice 
differentiable on  $\R-2\pi\Z$. Assume moreover the following
\begin{itemize}
\item [$(i)$] $\sup_{n\in \N} \int_{-\pi}^\pi  \frac{|\hat \sigma_n(\theta)|}{|\theta|} d\theta<\infty$;
\item [$(ii)$] $\sup_{n\in \N}\int_{-\pi}^\pi  |\hat \sigma_n'(\theta)| d\theta<\infty$;
 \item [$(iii)$] $\sup_{n\in \N} \int_{-\pi}^\pi |\theta|\, 
 |\hat \sigma_n ''(\theta)|d\theta <\infty\,$ .
\end{itemize}
Then, $\sup_{n\in \N} \|\sigma_n\|_{\ell^1(\Z) }<\infty$.
\end{prop}
 \noindent {\bf Proof.} We first notice that, by $(i)$,
 $$
\sup_{n\in \N}|\sigma_n (0)|\le \sup_{n\in \N}  \int_{-\pi}^\pi  |\hat \sigma_n(\theta)| d\theta<\infty
 $$Let $k\in \Z-\{0\}$.  We have
 \begin{gather*}
 \sigma_n(k)=\int_{-\pi}^\pi \hat\sigma_n(\theta) {\rm e}^{-ik\theta} 
 d\theta \\= \int_{-\pi/|k|}^{\pi/|k|} \hat\sigma_n(\theta) {\rm e}^{-ik\theta} d\theta +\int_{[-\pi,\pi]- [-\pi/|k|,\pi/|k|]}  \hat\sigma_n(\theta) 
 {\rm e}^{-ik\theta}\, d\theta \, .
 \end{gather*}
 Integrating by part and  using that $\hat \sigma_n$ is $2\pi$-periodic, we have
 \begin{gather*}
 \int_{[-\pi,\pi]- [-\pi/|k|,\pi/|k|]}  \hat\sigma_n(\theta) 
{\rm e}^{-ik\theta} \, d\theta 
 =-\int_{[-\pi,\pi]- [-\pi/|k|,\pi/|k|]}  \hat\sigma_n'(\theta) 
 \frac{{\rm e}^{-ik\theta}}{-ik} \, d\theta + \frac{\sigma_n(-\pi/|k|)
 -\sigma_n (\pi/|k|)}{-ik} \, ,
 \end{gather*}
and 
\begin{gather*}
\int_{[-\pi,\pi]- [-\pi/|k|,\pi/|k|]}  \hat\sigma_n'(\theta) 
 \frac{{\rm e}^{-ik\theta}}{-ik} \, d\theta=-\int_{[-\pi,\pi]- [-\pi/|k|,\pi/|k|]}  \hat\sigma_n''(\theta) 
 \frac{{\rm e}^{-ik\theta}}{-k^2} \, d\theta +\frac{\sigma_n'(-\pi/|k|)
 -\sigma_n' (\pi/|k|)}{-k^2} \,
\end{gather*}
Now,
\begin{gather*}
\sum_{|k|\ge 1} |\int_{-\pi/|k|}^{\pi/|k|} \hat\sigma_n(\theta) 
 {\rm e}^{-ik\theta}\, d\theta| \le 
 \int_{-\pi}^\pi \, |\hat \sigma_n(\theta)|\sum_{1\le |k|\le \pi/|\theta|}
 1 \, d\theta \le  2\pi \int_{-\pi}^\pi  \frac{|\hat \sigma_n(\theta)|
 }{|\theta|} d\theta
 \, ,
\end{gather*}
and 
$$
\sum_{|k|\ge 1}|\int_{[-\pi,\pi]- [-\pi/|k|,\pi/|k|]}  \hat\sigma_n''(\theta) 
 \frac{{\rm e}^{-ik\theta}}{-k^2} \, d\theta| \le 
 \int_{-\pi}^\pi |\hat \sigma_n ''(\theta)|\sum_{|k|\ge \pi/|\theta|} 
 \frac1{k^2}\le C  \int_{-\pi}^\pi |\theta|\, 
 |\hat \sigma_n ''(\theta)|d\theta \, .
$$

Hence, it remains to show that $\sup_{n\in \N} \sum_{|k|\ge 1} 
\frac{|\hat \sigma_n(\pi/k)|}{|k|}
 <\infty$ and  $\sup_{n\in \N}\sum_{|k|\ge 1} \frac{|\sigma_n '
 (\pi/k)| }{k^2}<\infty$.
 
 \smallskip
 
 Let $f_n(\theta):=\theta\hat \sigma_ n(\theta)$, for every $\theta
 \in \R-2\pi\Z$. Then $f_n$ is differentiable on $\R-2\pi\Z$ and, 
 by $(i)$ and $(ii)$, $\hat \sigma_n'\in L^1([0,2\pi])$, $f_n'\in L^1([0,2\pi])$. Hence, $\hat \sigma_n$ and $f_n$ can be continuously extended to $\R$ with $f_n(0)=0$. Then, for every $k\ge 1$, 
 \begin{gather*}
 \frac\pi{k}|\hat \sigma_n(\pi/k)|=\big|\int_0^{\pi/k}f_n'(\theta) d\theta \big|
 \le \int_0^{\pi/k}|\hat \sigma_n(\theta)|d\theta + \int_0^{\pi/k} \theta 
 |\hat \sigma_n'(\theta)|d\theta\, .
 \end{gather*}
 Dealing similarly with $k\le -1$ we infer that 
 \begin{gather*}
 \sum_{|k|\ge 1} 
\frac{|\hat \sigma_n(\pi/k)|}{|k|}\le \sum_{k\ge 1} \big( 
\int_{-\pi/k}^{\pi/k}|\hat \sigma_n(\theta)|d\theta + \int_{-\pi/k}^{\pi/k} \theta 
 |\hat \sigma_n'(\theta)|d\theta\big)\\ \le 
 \pi \int_{-\pi}^\pi  \frac{|\hat \sigma_n(\theta)|}{|\theta|} d\theta 
 +\pi\int_{-\pi}^\pi  |\hat \sigma_n'(\theta)| d\theta\, ,
 \end{gather*}
 which is bounded uniformly with respect to $n$. 
 
 \smallskip
 
 Proceeding as above with $g_n(\theta):=\theta^2\hat\sigma_n'(\theta)$ in 
 place of 
 $f_n(\theta)$ we see that, by $(ii)$ and $(iii)$, $\sup_{n\in \N}\sum_{|k|\ge 1} \frac{|\sigma_n '
 (\pi/k)| }{k^2}<\infty$. \hfill $\square$

 \medskip
 
 Let $\mu$ be a probability measure on $\Z$. We say that $\mu$ satisfies 
 hypothesis \Ht\,  if it satisfies hypothesis \H\, with a function $\psi$ such that there exists $D>0$ such that for every $\theta\in 
(0,\pi]$, 
\begin{equation}\label{psi-cond}
 \psi(\theta)\le D \theta \psi'(\theta)\, .
\end{equation}

\smallskip

\begin{prop}\label{rittprop}
Let $\mu$ be a probability measure on $\Z$ satisfying hypothesis \Ht. 
Then $(\sigma_n)_{n\in \N}:=(n(\mu^{*n}-\mu^{*(n+1)}))_{n\in \N}$ 
satisfies to items $(i)$, $(ii)$ and $(iii)$ of Proposition \ref{genprop2}. 
In particular, $\sup_{n\in \N}n\|\mu^{*n}-\mu^{*(n+1)}\|_{\ell^1(\Z)}
<\infty$, i.e. $\mu$ is Ritt.
\end{prop}
\noindent {\bf Proof.} By Proposition \ref{genprop} we already know 
that  $(iii)$ holds. It follows from the proof of Proposition \ref{genprop} 
and from \eqref{psi-cond} that there exist $C,c>0$ such that for every $\theta\in (0,\pi]$, 
\begin{gather*}
|\hat\sigma_n(\theta)|/\theta \le C n {\rm e}^{-cn\psi(\theta)} \psi(\theta)
/\theta\le CD n {\rm e}^{-cn\psi(\theta)} \psi'(\theta)\\
|\hat\sigma_n'(\theta)|\le C n{\rm e}^{-cn\psi(\theta)}\psi'(\theta)(n  
\psi(\theta)+ 1)\, .
\end{gather*}
then, we conclude as in the proof of Proposition \ref{genprop}. 
\hfill $\square$

\medskip

We now provide a sufficient condition on sequence of finite 
signed measure on $\Z$ to be bounded in $\ell^1(\Z)$, that will be needed 
in the sequel.

\begin{prop}\label{ritt-ext}
Let $(\sigma_n)_{n\in \N}$ be a sequence of finite signed measures 
on $\Z$ such that for every $n\in \N$, $\hat 
\sigma_n$ is continuously differentiable on $[-\pi,\pi]-\{0\}$. Assume that 
\begin{itemize}
\item [$(i)$] $\sup_{n\in \N} n\int_{-\pi}^\pi |\hat\sigma_n(\theta)
|d\theta <\infty$;
\item [$(ii)$] $\sup_{n\in \N} \int_{-\pi}^\pi \frac{|\hat\sigma_n'
(\theta)|^2}{n+1}
d\theta <\infty$.
\end{itemize}
Then, $\sup_{n\in \N}\|\sigma_n\|_{\ell^1(\Z)}<\infty$.
\end{prop}
\noindent {\bf Proof.} Let $n\ge 1$. Let $k\in \Z$. We have
\begin{gather}
\label{fourier1} \sigma_n(k)=\int_{-\pi}^\pi \hat \sigma_n(\theta) 
{\rm e}^{-ik\theta}d\theta\, ,
\end{gather}
and if $k\neq 0$, 
\begin{gather}
\label{fourier2} \sigma_n(k)=\int_{-\pi}^\pi \hat \sigma_n'(\theta) 
\frac{{\rm e}^{-ik\theta}}{ik}d\theta\, ,
\end{gather}
Using \eqref{fourier1}, we infer that 
$\sum_{0\le |k|\le n} |\sigma_n(k)|\le (2n+1) \int_{-\pi}^\pi |\hat\sigma_n(\theta)
|d\theta$. 
Using \eqref{fourier2}, Cauchy-Schwarz and Parseval, we infer that 
\begin{gather*}
(\sum_{|k|>n} |\sigma_n(k)|)^2 \le \Big(\int_{-\pi}^\pi |\hat \sigma_n'
(\theta) |^2 d\theta \Big)\sum_{|k|>n} \frac1{k^2}
\le \frac{C}{n}\int_{-\pi}^\pi |\hat \sigma_n'
(\theta) |^2 d\theta\, .
\end{gather*}
Then, we conclude thanks to $(i)$ and $(ii)$. \hfill $\square$

\subsection{Centered probability measures with a second moment}

It is known, see \cite{BC}, that a centered and strictly aperiodic probability 
measure $\mu$ on $\Z$ 
with a second moment satisfies \eqref{weakinemax0}. As an application 
of the previous subsections we add here that $\mu$ is moreover Ritt and satisfies \eqref{weakinemax}. Indeed, we shall prove Theorem \ref{centered}.

By Proposition \ref{genprop} and Proposition \ref{rittprop}, it suffices to prove that a centered and strictly 
aperiodic probability 
measure $\mu$ with a second moment satisfies condition \Ht \, for some 
function $\psi$. 

We shall take $\psi(\theta)=\theta^2$, for every $\theta\in [-\pi,\pi]$. 
Then $\psi$ satisfies \eqref{psi-cond} hence we just have to prove that 
$\mu$ satisfies \H.

Since $\mu$ has a second moment and is centered, it is twice continuously 
differentiable on $[-\pi,\pi]$ and we have 
$$
\lim_{\theta\to 0,\theta\neq 0}(1-{\rm Re}\, \hat\mu(\theta))/\theta^2 = 
\hat\mu''(0)/2>0\, ,
$$
and 
$$
\lim_{\theta\to 0,\theta\neq 0}{\rm Im}\, \hat\mu(\theta))/\theta^2 = 
0\, .
$$
It follows that item $(i)$ of hypothesis \H\, is satisfied for $\theta$ 
close enough to 0. Then, taking $c$ smaller if necessary, it holds on $(0,\pi]$ 
by strict aperiodicity. 

Using again that $\mu$ has a second moment and is centered we see that 
for every $\theta\in [-\pi,\pi]$, $|\hat \mu'(\theta)|\le \|\hat 
\mu''\|_\infty |\theta|$. Hence items $(ii)$ and $(iii)$ of hypothesis 
\H \, hold. Similarly, item $(iv)$ holds.

\section{ Probability measures  without first moment }

\bigskip

In this section as well as in sections 4 and 5, we shall consider probability measures $\mu$ on $\Z$ such that 
$(\mu(n))_{n\in \N}$ is \emph{completely monotone 
sequence}. Let us recall some definition and facts. 

\begin{defn}
Let $\Delta$ be the operator defined for every sequence $(t_n)_{n\in \N}$ 
of real numbers, by  $(\Delta t_n)_{n\in \N}=(t_n-t_{n+1})_{n\in \N}$. 
We say that a sequence $(t_n)_{n\ge 0}$ is completely monotone if 
for every $m\ge 0$ (with the convention $\Delta^0=Id$, 
$(\Delta ^mt_n)_{n\ge 0}$ is non-negative.
\end{defn}

\begin{defn}
We say that an infinitely differentiable function $f \,: \, 
[s,+\infty)\to [0,+\infty)$  is completely monotone, if 
for every $m\ge 0$, $(-1)^mf^{(m)}\ge 0$. 
\end{defn}

The following characterization of completely monotone sequences is due to Hausdorff 
and may be found in Widder \cite{Widder}, p.108.
\begin{prop}[Hausdorff]
A sequence $(\mu_n)_{n\in \N}$ is completely monotone if and only if there exists 
a finite positive measure $\nu$ on $[0,1]$, such that $\mu_n=\int_0^1
t^n \nu(dt$ for every $n\in \N$.
\end{prop}

A way to generate completely monotone sequences is the following, see 
\cite{Widder}, Theorem 11d, p. 158. 

\begin{prop}\label{widder}
Let $f$ be a completely monotone function. Then $(f(n+1))_{n\in \N}$ is a completely monotone sequence. 
\end{prop}

\medskip

\begin{defn}
We say that a probability measure $\mu$ on $\Z$ is CM if it is supported on $\N$ and if there exists 
 a finite (positive) measure  $\nu$ on $[0,1]$, such that 
\begin{equation}\label{proba}
\int_0^1 \frac{\nu(dt)}{1-t}=1\,. \footnote{All along the paper (for esthetical reasons) we shall adopt the convention $\int_a^b \varphi d\nu =\int_{[a,b]}\varphi d\nu$. Hence, for non-negative $\varphi$ 
we will have $\int_a^c \varphi d\nu \le \int_a^b \varphi d\nu+ \int_b^c \varphi d\nu$ with equality if $\nu(\{b\})=0$}
\end{equation}
and 
\begin{equation}\label{defi}
\mu(n)=\int_0^1 t^n \nu(dt)\qquad \forall n\in \N\,  .
\end{equation}
To emphasize the measure $\nu$ we shall say that $\mu$ is a CM probability measure on $\Z$ with representative measure $\nu$.
\end{defn}

\medskip

Notice that for $\mu$ as above, $\mu(n)>0$ for every $n\in \N$, hence 
$\mu$ is strictly aperiodic.

\subsection{Characterization of the BAR property} 

We first give an equivalent formulation of the BAR property that 
will be more convenient in the sequel. 

\begin{defn}
We say that a subset of ${\mathbb C}$ is a Stolz region if it is the convex hull 
of $1$ and a circle centered at $0$, with radius $0<r<1$.
\end{defn}

It is known that $\mu$ is strictly aperiodic and has BAR if and only if the range of $\hat \mu$ is 
included in a Stolz region.

\medskip

If $\mu$ is strictly aperiodic, for every $\varepsilon \in (0,\pi)$, 
$\hat \mu([\varepsilon, 2\pi-\varepsilon])$ is included in a disk centered at 0 with radius strictly smaller than 1. Hence, a strictly aperiodic 
$\mu$ has BAR if and only if 
\begin{equation}\label{sector}
\sup_{\theta\in (0,2\pi)}\frac{|{\rm Im}\, (\hat \mu(\theta))|}
{1-{\rm Re}\, (\hat \mu(\theta))|}
<\infty\, .
\end{equation}

%\medskip

%We wish to provide here a large class of   probabilities on $\N$ with bounded angular ratio. 

%Let us first notice that the set of probabilities 
%that have BAR is closed by convex combination and by convolution.

\medskip

We shall consider the following  condition on $\nu$: there 
exists $L>0$ such that for every $x\in [0,1)$,

\begin{equation}\label{condBAR}
\int_0^x \frac{t}{(1-t)^2}\,\nu(dt) \le  \frac{L}{1-x}
 \int_x^1 \frac{t}{1-t}\, \nu(dt)\,  .
\end{equation}
Notice that this condition implies that $\int_0^1\frac{\nu(dt}{(1-t)^2}=+\infty$ or, equivalently, that $\sum_{n\in \N} na_n=+\infty$, i.e. $\mu$ does not 
have first moment.

\begin{prop}\label{prop}
Let $\mu$ be a CM probability measure on $\Z$ with representative measure $\nu$. Then, $\mu$ has BAR if and only if there exists $L>0$ 
such that $\nu$ satisfies \eqref{condBAR}. Moreover, then
\begin{equation}\label{estimate}
\frac{1- {\rm Re }\, \hat\mu(\theta)}{ |\theta|}\underset{\theta\to 0}\longrightarrow +\infty
\, .
\end{equation}
\end{prop}

We deduce the following corollary, in the spirit of Theorem 4.1 of Dungey \cite{Dungey11}.

\begin{cor}\label{cor}
Let $\mu$ be a CM probability measure on $\Z$ with representative measure $\nu$ satisfying \eqref{condBAR} for 
some $L>0$.  Let $\tau$ be a probability measure on $\Z$ 
such that there exists $c>0$ such that 
$$
\sum_{n\in \Z} n|\tau(n)-a\mu(n)|<\infty\, .
$$
Then, $\tau$ has BAR.
\end{cor}
%\noindent{\bf Remark.} The condition \eqref{control} implies that 
 
 Throughout the paper we will make use of the following easy inequalities.
 \begin{gather*}
 |\sin \theta|\le |\theta| \quad , \quad 1-\cos \theta \le \frac{\theta^2}2\qquad \forall \theta\in \R\, ,\\
 |\sin \theta |\ge 2|\theta|/\pi \quad , \quad 1-\cos\theta \ge \frac{\theta^2}{4}\qquad \forall \theta\in [-1,1]\, .
 \end{gather*}

\noindent {\bf Proof of Proposition \ref{prop}.} Assume first that $\nu$ 
satisfies \eqref{condBAR}. Since $\nu$ is not null, the support of $\mu$ 
is $\N$ and $\mu$ is strictly aperiodic.

\medskip

Hence, we just have to prove that there exists 
$K>0$, such that 
\begin{equation}\label{seccond}
|{\rm Im} \, \hat \mu(\theta)|\le K (1-{\rm Re } \, \hat \mu(\theta)) \qquad \forall \theta\in[-\pi,\pi]\, .
\end{equation}
%i.e. that $\hat \mu([-\pi,\pi])$ is included in some sector. 

\medskip

We have, for every $\theta\in[-\pi,\pi]$, $\hat \mu(\theta)=\int_0^1 \frac{\nu(dt)}{1-t{\rm e}^{i\theta}}\, $. 

\medskip

Notice that $|1-t{\rm e}^{i\theta}|^2= 1+t^2-2t\cos \theta=
(1-t)^2+2t(1-\cos \theta)$ and that 

\begin{gather*}{\rm Re}\, \big(\frac{1}{1-t}- \frac{
(1-t{\rm e}^{-i\theta})}{|1-t{\rm e}^{i\theta}|^2}\big)= \frac{(1-t)^2+2t(1-\cos \theta) -(1-t)(1-t\cos \theta)}{
(1-t)|1-t{\rm e}^{i\theta}|^2}\\=\frac{t(1-\cos \theta)}{
(1-t)|1-t{\rm e}^{i\theta}|^2}\, .
\end{gather*}
 Hence, using \eqref{proba}, we have 
\begin{gather}\label{idremu}
1-{\rm Re } \, \hat \mu(\theta)= 
\int_0^1 \frac{t(1-\cos\theta)}{(1-t)((1-t)^2+2t(1-\cos \theta))}\,\nu(dt)\,  .
\end{gather}

Moreover, 
\begin{equation}\label{idimmu}
{\rm Im}\, \hat \mu(\theta) = \int_0^1 \frac{t\sin \theta}{|1-t{\rm e}^{i\theta}|^2}\, \nu(dt)= 
\int_0^1 \frac{t\sin \theta}{(1-t)^2+2t(1-\cos \theta)}\,\nu(dt) \, .
\end{equation}

Since, $\hat \mu$ is continuous and $1-{\rm Re}\, \hat \mu$ vanishes 
only at 0, on $[-\pi,\pi]$, it is enough to prove 
\eqref{seccond} for $\theta\in [-1/2,1/2]$. Moreover, \eqref{seccond} 
is clear for $\theta=0$. So, let $\theta\in
 [-1/2,1/2]-\{0\}$. 
 
 \medskip
 
 Let us first estimate $1-{\rm Re} \ \hat \mu(\theta)$.  Using that $(1-t)^2+
2t(1-\cos \theta)\le (1-t)^2+ \theta^2\le 2\max((1-t)^2,\theta^2)$, we obtain 
\begin{gather}\label{lowremu}
 1-{\rm Re}\, \hat \mu(\theta)\ge \frac12\int_0^{1-|\theta|} 
 \frac{t(1-\cos \theta)}{(1-t)^3}\nu(dt)\, +\, \frac18\int_{1-|\theta|}^1
 \frac{t}{1-t} \nu(dt)
\end{gather}

Now, we estimate ${\rm Im }\, \hat \mu$. We have, 
%using that for $t\ge 1-|\theta|\ge 1/2$, $(1-t)^2+2t(1-\cos \theta)\ge 1-\cos \theta$
\begin{gather*}
\int_{1-|\theta|}^1 \frac{t|\sin \theta|}{(1-t)^2+2t(1-\cos \theta)}\,\nu(dt) \le \int_{1-|\theta|}^1 \frac{t\theta^2}{(1-t)((1-t)^2+2t(1-\cos \theta))}\, \nu(dt)\\ \le 4 \int_{1-|\theta|}^1\frac{t(1-\cos \theta)}{(1-t)((1-t)^2+
2t(1-\cos \theta))}\,  \nu(dt)\le 4(1-{\rm Re}\, \hat\mu(\theta))\, .
\end{gather*}

Now, using our assumption on $\nu$ and \eqref{lowremu}, we obtain
\begin{gather*}
\int_0^{1-|\theta|} \frac{t|\sin \theta|}{(1-t)^2+2t(1-\cos \theta)}\,\nu(dt)
\le \int_0^{1-|\theta|} \frac{t|\sin \theta|}{(1-t)^2}\,\nu(dt)\\
\le \frac{L}{|\theta|} \int_{1-|\theta|}^1 \frac{t}{
1-t}\,\nu(dt)\le 8L (1-{\rm Re}\, \hat \mu(\theta))\, .
\end{gather*}
and we see that \eqref{seccond} holds.

\medskip

Let us prove the converse. Assume  that \eqref{seccond} holds.

Let $S\ge 1$ be fixed for the moment. Let $\theta\in [-1/2S,1/2S]-\{0\}$.

Using that $|1-t{\rm e}^{i\theta}|^2\le (1+1/S^2)(1-t)^2$, whenever 
$0\le t\le 1-S|\theta|$, we see that
\begin{equation}\label{min}
\int_0^{1-S|\theta|} \frac{\nu(dt)}{(1-t)^2}\le \frac{1+1/S^2}{2|\theta|/\pi} 
|{\rm Im}\, \hat 
\mu(\theta)| \le C \frac{1+1/S^2}{2|\theta|/\pi} (1-{\rm Re }\, \hat 
\mu(\theta))\, .
\end{equation}

Now, we see that
$$
1-{\rm Re }\, \hat \mu(\theta)) \le \frac{1-\cos\theta}{S|\theta|} 
\int_0^{1-S|\theta|}\frac{t\nu(dt)}{(1-t)^2} +\int_{1-S|\theta|}^1 
\frac{t\nu(dt)}{1-t}\, .
$$

 Hence, taking $S$ large enough and using \eqref{min}, we infer that 
 there exists $D>0$ such that
 $$
 \int_0^{1-S|\theta|} \frac{\nu(dt)}{(1-t)^2}\le \frac{D}{|\theta|}
 \int_{1-S|\theta|}^1 
\frac{t\nu(dt)}{1-t}\, ,
 $$
 which prove that \eqref{condBAR} holds \hfill $\square$

\medskip

It remains to prove \eqref{estimate}. Using \eqref{condBAR}, we see that
\begin{gather*}
\frac{{ 1-\rm Re}\,\hat \mu (\theta)}{|\theta|}\ge \frac{1-\cos\theta}
{2|\theta|^3} \int_{1-|\theta|} ^1
\frac{t}{(1-t)}\, \nu(dt) \ge \frac{1-\cos\theta}
{2|\theta|^2} \int_0^{1-|\theta|} 
\frac{t}{(1-t)^2}\, \nu(dt)
\underset{\theta\to 0}\longrightarrow +\infty\, ,
\end{gather*} 
hence the result.
\medskip

 \hfill $\square$

\medskip

\noindent {\bf Proof of Corollary \ref{cor}.} By assumption and Proposition 
\ref{prop}, there exists 
$K>0$ such that $\sum_{n\ge 1}n|\tau(n)-a\mu(n)|\le K$ and for 
every $\theta\in [-\pi,\pi]$, $|{\rm Im}\, (\hat\mu(\theta))|\le 
K(1-{\rm Re} \, (\hat \mu(\theta ))$.
% and $|\theta|\le K(1-{\rm Re} \, (\hat \mu(\theta ))$.

\medskip

Let us  prove that $\tau$ is  strictly aperiodic. If $\tau$ were not 
strictly aperiodic, there would exists $\ell\ge 2$ and $0\le k\le \ell-1$, such that the support of $\tau$ would be contained in $k+\ell \Z$. 
In particular $\tau(k+1+\ell m)=0$ for every $m\in \Z$. Hence, 
$\sum_{m\in \Z}|m| \mu(k+1+\ell m)<\infty$ and (using that $(\mu(n))_{n\ge 1}$ is non increasing) $\mu$ must have a first 
moment, contradicting \eqref{condBAR} (see the remark after 
\eqref{condBAR}).

\medskip

We first prove that there exists $C>0$ such that 
for every $\theta\in [-\pi,\pi]$, 
\begin{equation}\label{esti}
{\rm Re}\, (1-\hat \tau(\theta))\ge C|\theta|\, .
\end{equation} 
Since $\tau$ is strictly aperiodic, it is enough to prove the result for small 
enough $\theta$'s. By Proposition \ref{prop}, there exists $\delta\in (0,\pi)$, 
such that for every $\theta\in [-\delta,\delta]$, $|\theta|\le a(1-{\rm Re} \, (\hat \mu(\theta ))/2K$. Then, using that $1-\cos u\le |u|$ for every $u\in \R$,
\begin{gather*}
|\theta|\le (1-{\rm Re} \, (\hat \tau(\theta ))/2K+\frac1{2K} 
\sum_{n\ge 1}|\tau(n)-a\mu(n)|(1-\cos(n\theta))\\ \le (1-{\rm Re} \, (\hat \tau(\theta ))/2K+|\theta|/2 \, ,
\end{gather*}
 and \eqref{esti} follows.

\medskip Let $\theta\in [-\pi,\pi]$. We have, using that 
$|\sin u|\le u$  for every $u\in \R$,
\begin{gather*}
|{\rm Im}\, (\hat\tau(\theta))|\le a|{\rm Im}\, (\hat\mu(\theta))| 
+\sum_{n\ge 1} |\tau(n)-a\mu(n)|\, |\sin(n\theta) |
\le aK (1-{\rm Re} \, (\hat \mu(\theta )) + K|\theta|
\\ \le K (1-{\rm Re} \, (\hat \tau(\theta )) + \sum_{n\ge 1} 
|\tau(n)-a\mu(n)|\, |(1-\cos (n\theta)| +K|\theta| \\
\le K (1-{\rm Re} \, (\hat \tau(\theta )) +2K|\theta|\le 
K (1+2C)(1-{\rm Re} \, (\hat \tau(\theta )) \,,
\end{gather*}
and the corollary is proved. \hfill $\square$

\medskip

\medskip

From a pratical point of view it is better to have a condition on $(\mu(n))_{n\in \Z}$. Indeed we may consider completely monotone 
sequences given thanks to Proposition \ref{widder}, in which 
case, we do not know $\nu$. 

\medskip

\begin{prop}\label{propcarac}
Let $\mu$ be a CM probability measure on $\Z$ with representative measure $\nu$. Then,  $\nu$  satisfies \eqref{condBAR}
 if and only if 
  there exists $D>0$, such that for every $n\ge 1$, 
 \begin{equation}\label{condBAR2}
\sum_{k=1}^n k\mu(k) \le Dn\sum_{k\ge n} \mu(k) \, . 
 \end{equation}
\end{prop}
\noindent {\bf Proof.} 
Assume \eqref{condBAR}. Let $n\ge 1$. We have 
\begin{gather*}
\sum_{k=1}^n k\mu(k) \le \int_0^{1-1/n} \frac{t}{(1-t)^2}\nu (dt) 
+ n \int_{1-1/n}^1 \frac{t}{1-t}\nu (dt) \le (1+L)n 
\int_{1-1/n}^1 \frac{t}{1-t}\nu (dt)\, .
\end{gather*}
Using that $\sum_{k\ge n} \mu(k)=\int_0^1 \frac{t^n}{1-t}\nu(dt)$ and that 
$(1-1/n)^n\underset{n\to \infty}\longrightarrow {\rm e}^{-1}$, we see that 
\eqref{condBAR2} holds. 

\medskip 

Assume now that \eqref{condBAR2} holds. 

Let $A\ge 1$ be a positive integer fixed for the moment. Let $n\ge 2$. 
\medskip

Let $1\le m\le n-1$ be an integer and let $t\in [1-1/m,1-1/(m+1)]$. Using that 
the sequence $((1-1/k)^{k-1})_{k\ge 1}$ decreases to $1/{\rm e}$, we obtain that 
(with the convention $0^0=1$)
$$
\sum_{k=1}^{An}kt^k \ge t\sum_{k=0}^{m-1}(k+1)(1-1/m)^{m-1}\ge \frac{tm(m+1)}{2
{\rm e}}\ge \frac{t}{{\rm e}(1-t)^2}\, .
$$

Hence, 
$$
\int_0^{1-1/n} \frac{t}{(1-t)^2}\nu(dt)\le {\rm e} \sum_{k=1}^{An} k\mu(k)
\le {\rm e}DAn\sum_{k\ge An} \mu(k) = {\rm e}DAn\int_0^1 \sum_{k\ge An}
t^k\nu(dt)\, .
$$

Now notice that for $t\in [0,1-1/n]$,  $$\sum_{k\ge An} t^{k-1} \le 
\frac1{A^2n^2}\sum_{k\ge An} k(k+1)t^{k-1}\le \frac1{A^2n^2(1-t)^3} 
\le \frac1{A^2n(1-t)^2}\, ,$$
and that for $t\in [1-1/n,1]$, $\sum_{k\ge An} t^{k-1}\le 1/(1-t)$.

Hence, taking $A$ large enough we infer that  \eqref{condBAR} holds. \hfill $\square$

\subsection{Hypothesis ${\bf (H)}$ for CM probability measures}

We shall prove  that the conditions imposed in 
the previous subsection guarantee hypothesis ${\bf (H)}$. 

\medskip

\begin{prop}\label{theocomp}
Let $\mu$ be a CM probability  measure on $\Z$ satisfying \eqref{condBAR2}. Then, $\mu$ satisfies hypothesis \H. 
%\label{weakinemax}
%\begin{equation*}
%.\sup_{\p >0} \p \#\{k\in \Z\,:\, \sup_{n\ge 1} n^m|\mu^{*n}*(\delta_0-\mu)^{*m}*f(k)|\ge \p\}
%\le C_m\|f\|_{\ell^1}\, .
%\end{equation*}
\end{prop}
\noindent {\bf Proof.} 
%Since we already know that $\mu$ is Ritt, by Proposition \ref{genprop}, the proof will be completed if we can prove that 
%$\mu$ satisfies hypothesis \H.
%\medskip
To check the conditions $(i)-(iv)$ of hypothesis \H\, with a suitable function 
$\psi$ we must first estimate $\hat \mu$ and its derivatives. 
 
%As we will see those bounds concerning the real part of $\hat \mu$ do not make use of the fact that 
%$(\mu(n))_{n\in \N}$ satisfies  \eqref{condBAR2}. 

\medskip

Let us first compute the derivatives of $\hat \mu$. Recall that for every 
$\theta\in [-\pi,\pi]$,

\begin{gather*}
 1-{\rm Re}\, \hat\mu(\theta) =\int_0^1 \frac{t(1-\cos \theta)}
 {(1-t)\big((1-t)^2+2t(1-\cos \theta)\big)} \, \nu(dt) \\
 = \frac12\int_0^1 \frac{\nu(dt)}{1-t} \, -\, \frac12 \int_0^1 
 \frac{(1-t)}{(1-t)^2+2t(1-\cos \theta)}\, \nu(dt)\, ,
 \end{gather*}
 and 
 \begin{gather*}
{\rm Im}\, \hat \mu(\theta) = \int_0^1 \frac{t\sin \theta}{|1-t{\rm e}^{i\theta}|^2}\, \nu(dt)= \int_0^1 \frac{t\sin \theta}{
(1-t)^2+2t(1-\cos \theta)}\, \nu(dt)\, .
\end{gather*}

 Hence, for every $\theta\in [-\pi,\pi]-\{0\}$,
  \begin{gather}
\label{remu'} {\rm Re}\, \hat\mu'(\theta) = -\sin \theta \int_0^1 
 \frac{t(1-t)}{\big((1-t)^2+2t(1-\cos \theta)\big)^2}\, \nu(dt)\, ,
 \end{gather}
 \begin{gather}
\label{remu''} {\rm Re}\, \hat\mu''(\theta) = -\cos \theta \int_0^1 
 \frac{t(1-t)}{\big((1-t)^2+2t(1-\cos \theta)\big)^2}\, \nu(dt)\,
 \\ \nonumber+2\sin^2 \theta \int_0^1 
 \frac{t(1-t)}{\big((1-t)^2+2t(1-\cos \theta)\big)^3}\, \nu(dt)
 \, 
 \end{gather}

\begin{gather}
\label{immu'}{\rm Im}\, \hat \mu'(\theta) = \int_0^1 \frac{t\cos \theta}{
(1-t)^2+2t(1-\cos \theta)}\, \nu(dt)\qquad  \\ \nonumber \qquad \qquad -
\int_0^1 \frac{2t^2\sin^2 \theta}{\big(
(1-t)^2+2t(1-\cos \theta)\big)^2}\, \nu(dt)\, ,
\end{gather}

\begin{gather}
\label{immu''}{\rm Im}\, \hat \mu''(\theta) = \int_0^1 \frac{-t\sin \theta}{
(1-t)^2+2t(1-\cos \theta)}\, \nu(dt)\qquad  \qquad \\
\nonumber \qquad \qquad -
\int_0^1 \frac{8t^2\sin \theta\cos \theta }{\big(
(1-t)^2+2t(1-\cos \theta)\big)^2}\, \nu(dt)\,  +\int_0^1 \frac{4t^3\sin^3 \theta }{\big(
(1-t)^2+2t(1-\cos \theta)\big)^3}\, \nu(dt)\,.
\end{gather}

Define, for $\theta\in [-\pi,\pi]$, 
\begin{equation}\label{psi-def}\psi(\theta)=
\int_0^1\frac{t|\theta|}{(1-t)(1-t+t|\theta|)}\,\nu(dt)
=1-\int_0^1\frac{\nu(dt)}{(1-t+t|\theta|)}\, . 
\end{equation}
 Hence, for every $\theta\in (0,\pi]$,

\begin{gather}\label{psi'-def}
\psi'(\theta)=\int_0^1 \frac{t\nu(dt)}{
(1-t+t\theta)^2}\, .
\end{gather}

Notice that, for every $\theta\in (0,1/2]$, 

\begin{gather}
\label{psi} \frac{\theta}2 \int_0^{1-\theta}
\frac{t\nu(dt)}{(1-t)^2}\,+\,\frac{1}{2} \int_{1-\theta}^1
\frac{t\nu(dt)}{1-t}\le \psi(\theta)\le 
\theta \int_0^{1-\theta}
\frac{t\nu(dt)}{(1-t)^2}\,+\, \int_{1-\theta}^1
\frac{\nu(dt)}{1-t}\, ;\\
\label{psi'} \frac\theta4\int_0^{1-\theta} \frac{t\nu(dt)}{
(1-t)^2}\, +\, \frac{1}{4\theta}\int_{1-\theta}^1t\nu(dt)\le 
 \theta\psi'(\theta) \le \theta \int_0^{1-\theta}
\frac{t\nu(dt)}{(1-t)^2}\,+\,\frac{2}{|\theta|} \int_{1-\theta}^1
\nu(dt)\, .
\end{gather}

\begin{claim}\label{claim1}
There exists $C>0$, such that for every $\theta\in [0,\pi]$,
\begin{equation*}
1-{\rm Re}\, \hat\mu(\theta)\ge C \psi(\theta)\, .
\end{equation*}
\end{claim}
\noindent {\bf Proof.} It suffices to prove the claim 
for $\theta\in(0,1/2]$. Using \eqref{lowremu}, \eqref{condBAR} and 
\eqref{psi}  we have 
\begin{gather*}
1-{\rm Re}\, \hat\mu(\theta)\ge  
\int_{1-|\theta|}^1 \frac{t}{8(1-t)}\nu(dt)
\ge \frac{1}{8(L+2)}\psi(\theta)\, ,
\end{gather*}
and the claim follows.\hfill $\square$

\begin{claim}\label{claim2}
There exists $C>0$, such that for every $\theta\in (0,
\pi]$, $| \hat\mu'(\theta)|\le C \psi'(\theta)\, .$
\end{claim}
\noindent {\bf Proof.}  Again, we only consider the case when 
$\theta\in(0,1/2]$. We deal separately with the real and imaginary 
part of $\mu'$. We have, using \eqref{remu'} and \eqref{psi'}
\begin{equation*}
|{\rm Re }\, \hat\mu'(\theta)|\le 2\theta\int_0^{1-\theta}\frac{t\nu(dt)}
{(1-t)^3} \, +\, \frac2{\theta^3}\int_{1-\theta}^1t(1-t) \nu(dt)\, 
\le 8 \psi'(\theta)\, ,
\end{equation*}

Similarly, 
\begin{gather*}
|{\rm Im }\, \hat\mu'(\theta)|\le \int_0^{1-\theta} \Big(\frac{t}
{(1-t)^2}+\frac{t\theta^2}{(1-t)^4}\Big)\, \nu(dt) \, + \, 
\Big(\frac{1}{1-\cos \theta} + \frac{2\theta^2}{(1-\cos \theta)^2}\Big)
\int_{1-\theta}^1t\nu(dt)\\
\le C \psi'(\theta)\, .
\end{gather*}

\smallskip

\begin{claim}\label{claim3}
There exists $C>0$, such that for every $\theta\in (0,
\pi]$, $|\theta  \hat\mu'(\theta)|\le C \psi(\theta)\, .$
\end{claim}
\noindent {\bf Proof.} Combine Claim \ref{claim2} and \eqref{psi'}.
\hfill $\square$

\medskip

\begin{claim}\label{claim4}
There exists $C>0$, such that for every $\theta\in (0,
\pi]$, $| \theta\hat\mu''(\theta)|\le C \psi'(\theta)\, .$
\end{claim}
\noindent {\bf Proof.} We assume that $\theta\in (0, 1/2]$. By \eqref{remu''} and \eqref{psi'}, we have 
\begin{gather*}
|{\rm Re }\, \hat\mu''(\theta)|\le \int_0^{1-\theta} \Big(\frac{2t}{(1-t)^3}+ \frac{4t\theta^2}{(1-t)^5}\Big)\, \nu(dt)\, \\
+\, \frac{2}{(1-\cos\theta)^2} \int_{1-\theta}^1 \frac{t}{1-t}\nu(dt)
+ \frac{4\theta^2}{(1-\cos\theta)^3}\int_{1-\theta}^1 
t(1-t)  \nu(dt)\, \le C\psi'(\theta)/\theta.
\end{gather*}
Similar computations based on \eqref{remu'} and \eqref{psi'} yields

\begin{gather*}
|{\rm Im}\, \hat\mu''(\theta)|\le  C\psi'(\theta)/\theta\, .
\end{gather*}

Then, items $ii)$, $(iii)$ and $(iv)$ of Proposition \ref{genprop}  follows from the combination 
 of the  claims 2, 3 and 4. 
 
 \smallskip
 
 Let us prove item $(i)$  of Proposition \ref{genprop}.
Let $\theta\in [0,1/2]$. Recall that $|{\rm Im} \, \hat \mu(\theta)|\le 
C(1-{\rm Re }\, \hat\mu(\theta))$. Hence, 
\begin{gather*}
|\hat \mu(\theta)|^2 = |{\rm Im} \, \hat \mu(\theta)|^2
+ 1 -2 (1-{\rm Re }\, \hat\mu(\theta))+(1-{\rm Re }\, \hat\mu(\theta))^2\\
1-(1-{\rm Re }\, \hat\mu(\theta))(2-(C+1)(1-{\rm Re }\, \hat\mu(\theta)))
\end{gather*} 
Since $(1-{\rm Re }\, \hat\mu(\theta))\underset{\theta\to 0}\longrightarrow 
0$, using Claim 1, we infer that for $\theta$ small enough
$$
|\hat \mu(\theta)|^2\le 1-\tilde C \psi(\theta)\, .
$$
Hence $(i)$ holds for, say,  $\theta\in [0,\eta]$, with $\eta$ small enough. Since $\sup_{\theta\in [\eta ,\pi]}|\hat \mu(\theta)|<1$, 
wee see that $(i)$ holds for every $\theta\in[0,\pi]$, taking $c$ smaller if necessary. 

 \hfill $\square$
 
 \smallskip
 
% Using that the function $\psi$ above satisfes $\liminf_{\theta\to 0,\theta>0}
%\psi'(\theta)=+\infty$, it is not hard to derive the following "perturbation" 
%result (see Proposition \ref{proppertbis} for some arguments). 

\begin{cor}\label{cor-CM}
Let $\mu$ be a CM probability measure on $\Z$. Let  $\sigma$ be a probability measure on $\Z$ such that $\hat \sigma$ is twice continuously differentiable 
on $[-\pi,\pi]-\{0\}$ and such that  $\hat \sigma'$  and $\theta\mapsto 
\theta \hat \sigma ''(\theta)$ are bounded. Then, $\sigma*\mu$ satisfies 
hypothesis ${\bf (H)}$. Moreover, if $\mu$ satisfies 
hypothesis ${\bf (\tilde H)}$, so does $\sigma*\mu$.
\end{cor}
{\bf Remark.} The assumptions on $\nu$ holds, for instance, as soon as 
$\sum_{n\in \Z}n^2\sigma(n)<\infty$.

\noindent {\bf Proof.} Let $\psi$ be the function defined in 
\eqref{psi-def}. Since $\int_0^1 \frac{\nu(dt)}{(1-t)^2}=+\infty$, one easily infers from \eqref{psi'-def} that 
$\liminf_{\theta\to 0,\theta>0}
\psi'(\theta)=+\infty$. In particular, there exists 
$K>0$ such that for every $\theta\in (0,\pi]$, 
$\psi'(\theta)\ge K$ and, consequently, 
$\psi(\theta)\ge K\theta$. Then, the fact that $\sigma*\mu$ satisfies 
hypothesis ${\bf (H)}$, with the same function $\psi$ as $\mu$,  may be proved exactly as Proposition \ref{ext}. Since we use the same function 
$\psi$ for $\sigma*\mu$ and $\mu$, then $\sigma*\mu$ satisfies 
hypothesis ${\bf ( \tilde H)}$ as soon as $\mu$ does. \hfill $\square$

\begin{cor}
Let $\tau$ be a probability measure on $\Z$. Assume  that 
there exists a CM probability measure  $\mu$  and $c>0$ such that 
$\sum_{n\in \Z}n^2|\tau(n)-c\mu(n)|<\infty$. Then,
$\tau$ satisfies hypothesis ${\bf (H)}$. If moreover $\mu$ 
satisfies hypothesis ${\bf (\tilde H)}$, so does $\tau$.
\end{cor}
\noindent {\bf Remark.} It follows from the proof that 
we only need that $\hat \sigma$ be twice continuously differentiable 
on $[-\pi,\pi]-\{0\}$ and that $\hat \tau'$ and $\theta\mapsto \theta \hat \tau''(\theta)$ be bounded.

\noindent {\bf Proof.} Define a \emph{signed} measure by setting $\sigma:=
\tau-c\mu$. Then, $\hat \sigma$ is twice continuously differentiable 
on $[-\pi,\pi]$, $\hat \sigma (0)=1-c$ and there exists 
$C>0$ such that for every $\theta\in [0,\pi]$,  
$|\hat \sigma(\theta)-(1-c)|\le C \theta$. Then, the proof may be finished using the same arguments as in the proof of Corollary 
\ref{cor-CM}. \hfill $\square$

%This corollary provides an improvement  of Proposition \ref{theocomp} 
%in the case where $m=1$. We cannot do the same when $m\ge 1$ due to the fact that 
%we are not aware of any perturbation version of \ref{propritt}. 

\subsection{The Ritt property on $\ell^1(\Z)$}
In this section, we finish the proof of Theorem \ref{theoNC}. We first prove the Ritt property of CM probability measures, which corresponds to the 
case where $\sigma=\delta_0$. 

\medskip

Let $\mu$ be a probability measure on $\Z$.  Notice that the fact 
$\mu$ is Ritt is equivalent to the fact that

$$
\sup_{n\ge 1} n\|\pi_\mu^{n}-\pi_\mu^{n+1}\|_{\ell^1(\Z)}<\infty \, ,
$$
where $\pi_\mu$ stands for the operator of convolution by $\mu$. 

\smallskip

Let $\Gamma$ be the open unit disk in the complex plane. By Theorem 1.5 of Dungey, $\mu$ is Ritt if and only if 
the spectrum $\sigma(\pi_\mu)$ of $\pi_\mu$ is contained in $\Gamma\cup \{1\}$ and the semi-group $({\rm e}^{-t(I-\pi_\mu)})_{t\ge 0}$ 
is bounded analytic. The fact that $({\rm e}^{-t(I-\pi_\mu)})_{t\ge 0}$ 
is bounded analytic means that

$$
\sup_{t>0} \Big( \|{\rm e}^{-t(\delta_0-\mu)}\|_{\ell^1(\Z)} +t\|(I-T) 
{\rm e}^{-t(\delta_0-\pi_\mu)} \|_{\ell^1(\Z)} \Big)  <\infty \, .
$$

\noindent {\bf Remark.} Notice that Theorem 1.5 of Dungey is valid for 
probabilities supported on $\N$. 

\medskip

\begin{prop}\label{propritt}
Let $\mu$ be a CM probability measure on $\Z$ with representative measure $\nu$ satisfying \eqref{condBAR}. Then, $\mu$ is Ritt.
\end{prop} 

We already saw that $\nu$ satisfies \eqref{condBAR} if and only if $\mu$ has BAR. The fact that a CM probability measure on $\Z$ having BAR is Ritt has been proved very recently (see their Theorem 7.1) by Gomilko and Tomilov \cite{GT1} 
as a consequence of another very recent   result of their own \cite{GT}. 

The latter paper deals with subordination semi-groups hence is written 
in a continuous setting.

For reader's convenience we explain below how to derive Proposition \ref{propritt} from the work \cite{GT}.

\medskip

First of all, by Theorem 2.1 of Dungey \cite{Dungey11}, we have 
$\sigma(\pi_\mu)\subset \hat \mu([-\pi,\pi])\subset \Gamma \cup\{1\}$, 
where the latter inclusion follows from the fact that $\mu$ has BAR. 
Hence, Proposition \ref{propritt} will be proved if we can prove that 
$({\rm e}^{-t(I-\pi_\mu)})_{t\ge 0}$ 
is bounded analytic.

%As in the previous section, let us assume that $\mu_n=\int_0^1 
%t^n\nu(dt)$, for every $n\in \N$, for some finite positive measure $\nu$ such that 
%$\int_0^1 \frac{\nu(dt)}{1-t}=1$. 

\medskip

\begin{defn}
An infinitely differentiable  function $f\, :\, (0,+\infty)\to [0,+\infty)$ is called a Bernstein function if $f'$ is completely monotone. If $\lim_{x\to 0^+}
f(x)$ exists and if $f$ admits an holomorphic extension to 
$\{z\in \C \, :\, {\rm Im}\, z>0\}$, such that ${\rm Im}\,f(z)\ge 0$, then
$f$ is called \emph{complete Bernstein}.
\end{defn}

For every $x\ge 0$, define $\chi(x):= 1- \int_0^1 \frac{\nu(dt)}{1-t+tx}
=\int_0^1 \frac{\nu(dt)}{1-t}-\int_0^1 \frac{\nu(dt)}{1-t+tx}$. 
Then $\chi$ is non-decreasing, with $\chi(0)=0$, hence it is 
non-negative. It is not hard to see that it is infinitely differentiable 
and that $\chi'$ is completely monotone, hence $\chi$ is a 
 Bernstein function and one can easily see that it is actually a 
 complete Bernstein function.

\medskip

Since $\chi$ is Bernstein, it is well known  (see e.g. Theorem 
1.2.4 of \cite{FJS}) that there exists a convolution semi-group 
$(\sigma_t)_{t\ge 0}$ (of probability measures on $[0,\infty)$), such that for every $x\ge 0$, and every $t\ge 0$, 
$$
\int_0^\infty {\rm e}^{-xy}\sigma_t(dy)={\rm e}^{-t\chi(x)}\, .
$$

\smallskip

Following Dungey \cite[p. 1734]{Dungey11}, we consider the Poisson semi-group 
$(P_s)_{s\ge 0}$ acting by convolution on $\ell^1(\N)$, and defined by 
$$
P_s:= {\rm e}^{-s(\delta_0-\delta_1)}= 
{\rm e}^{-s}\sum_{k\ge 0} \frac{s^k}{k!} \delta_k\qquad \forall s\ge 0\, .
$$

Consider now the associated subordinated semi-group $(Q_t)_{s\ge 0}$ defined by 
$$
Q_t:= \int_0^\infty P_s \, \sigma_t(ds)\qquad \forall t\ge 0\, .
$$

\smallskip

Let $t\ge 0$. Then, $Q_t$ is a probability measure on $\N$, whose generating function is given (on $[0,1]$) by 
$$
x\mapsto \int_0^\infty {\rm e}^{-s(1-x)}\sigma_t(ds)= {\rm e}^{-t\chi(1-x)}\, .
$$

Let $G_\mu$ denote the generating function of $\mu$, i.e. 
$$G_\mu(x)=\sum_{n\ge 0}\mu(n) x^n=\int_0^1\frac{\nu(dt}{1-tx}=1-\chi(1-x)\, ,
$$ for every $x\in [0,1]$. Then, for every $t\ge 
0$, the generating function of the probability 
${\rm e}^{-t(I-\mu)}= {\rm e}^{-t}\sum_{k\ge 0}\frac{t^k\mu^{*k}}{k!}$ is 
given by 
$$
{\rm e}^{-t}\sum_{k\ge 0}\frac{t^kG_\mu^{k}}{k!}={\rm e}^{-t(1-G_\mu)}\, .
$$

In particular, we see that the semi-groups $({\rm e}^{-t(I-\pi_\mu)})_{t\ge 0}$ 
and $(Q_t)_{t\ge 0}$ co\"\i ncide. Hence, to prove that $({\rm e}^{-t(I-\pi_\mu)})_{t\ge 0}$ 
is bounded analytic, it is enough to prove that any subordinated semi-group 
associated with $(\sigma_t)_{t \ge 0}$ is bounded analytic (see the introduction 
of \cite{GT} for more details). 
To prove the latter point, since $\chi$ is complete Bernstein, by 
Corollary 7.10 of \cite{GT}, it is enough to prove that 
$\chi$ sends the half-plane $\{z\in \C \, :\, {\rm Re}\, z\ge 0\}$ 
to a sector $\{ z\in \C \, :\, |{\rm Im} \, z|\le C{\rm Re}\, z\}$, for some $C>0$.

\medskip

 Let $z=a+ib$ such that $a\ge 0$ and $|z|^2=a^2+b^2\le 1/4$. 
We have, using \eqref{condBAR} 
\begin{gather*}
|{\rm Im}\, \chi(z)| = |b |\int_0^1 \frac{t}{(1-t+at)^2+t^2b^2}\nu(dt) 
\\ \le |b| \int_0^{1-|z|}  \frac{t}{(1-t)^2} \nu(dt) + \frac{4b}
{|z|^2} 
\int_{1-|z|}^1 t\nu (dt)
\le K\int_{1-|z|}^1 \frac{t}{1-t}\nu (dt)\, .
\end{gather*}
On the other hand, 
\begin{gather*}
{\rm Re}\, \chi(z)= \int_0^1 \frac{at(1-t)+|z|^2t^2}{(1-t)\big((
1-t+at)^2+t^2b^2\big)}\, \nu(dt)\ge \int_{1-|z|}^1
\frac{|z|^2 t^2}{5|z|^2(1-t)}\, \nu(dt) \\
\ge \frac1{10}   \int_{1-|z|}^1 \frac{t}{1-t}\nu (dt)\, .
\end{gather*}
This gives the desired bound when $|z| ^2\le 1/4$. 

\smallskip

Assume now that $|z|^2\ge 1/4$. In particular, we have $4|z|\ge 2$.  
Hence, 
\begin{gather*}
|{\rm Im}\, \chi(z)|  \le  \int_0^{(4|z|)^{-1}}  \frac{t|z|}{(1-t)^2} \nu(dt) + \frac{|z|}
{|z|^2} 
\int_{(4|z|)^{-1}}^1 t^{-1}\nu (dt)\\
\le \frac14 \int_0^{1/2} \frac{\nu(dt)}{(1-t)^2} + 4\int_0^1 \nu(dt)
<\infty \, .
\end{gather*}
Moreover, using that the integrand below is non decreasing with respect to 
$|z|$, we have
\begin{gather*}
{\rm Re}\, \chi(z)\ge  \int_0^1 \frac{|z|^2t^2}{2(1-t)\big((
1-t)^2+t^2|z|^2\big)}\, \nu(dt) \\
\ge \frac1{8}   \int_{0}^1 \frac{t^2}{(1-t)^2+t^2/4}\nu (dt)>0\, ,
\end{gather*}
which finishes the proof. \hfill $\square$

\begin{prop} \label{prop-conv}
Let $\mu$ be a CM probability measure on $\Z$ with representative measure $\nu$ satisfying \eqref{condBAR}. Let $\sigma$ be a probability measure 
on $\Z$ such that $\hat \sigma$ is continuously differentiable on 
$[-\pi,\pi]-\{0\}$ and such that $\hat\sigma '$ is bounded 
on $[-\pi,\pi]-\{0\}$. Then, 
\begin{equation}\label{int-ritt}
\sup_{n\in \N} n\|(\delta_0-\sigma)*\mu^{*n}\|_{\ell^1}<\infty\, .
\end{equation}
In particular, $\sigma*\mu$ is Ritt and for every $\alpha\in (0,1]$, 
$\alpha\mu +(1-\alpha )\sigma$ is Ritt. 
\end{prop}
\noindent {\bf Proof.} To prove \eqref{int-ritt}, we check that 
$(\sigma_n)_{n\in \N}:=(n(\delta_0-\sigma)*\mu^{*n})_{n\in \N}$ 
satisfies items $(i)$ and $(ii)$  of Proposition \ref{ritt-ext}. 

\smallskip

By assumption there  exists $L>0$ such that $|\hat\sigma'|\le L$ and it follows that $|1-\hat \sigma(\theta)|\le L|\theta|$ for every 
$\theta\in [-\pi,\pi]$.

\smallskip
Let $\psi$ be the function given in \eqref{psi-def}. Recall that there exists $K>0$ such that 
for every $\theta\in [-\pi,\pi]-\{0\}$, $\psi(\theta)\ge K\theta$ and 
$\psi'(\theta)\ge K$.  Hence, for every $n\in \N$ and every 
$\theta\in [-\pi,\pi]-\{0\}$,
$$
n|\hat \sigma_n(\theta)|\le 
\frac{Ln^2}{K^2} \psi(\theta)\psi'(\theta)(\hat \mu(\theta))^n\, .
$$
Hence, arguing as in the proof of Proposition \ref{rittprop}, 
we see that $(\sigma_n)_{n\in \N}$ satisfies item $(i)$ in Proposition 
\ref{ritt-ext}. 

\smallskip

For every $\theta\in [-\pi,\pi]-\{0\}$, we have 
\begin{equation*}
\hat\sigma_n'(\theta)=-n \hat\sigma' (\theta)\hat \mu^{n}(\theta)+ n^2 (1-\hat\sigma(\theta)) 
\hat\mu '(\theta)\hat \mu^{n-1}(\theta)\, .
\end{equation*}
Then, we infer that
\begin{equation*}
|\hat\sigma_n'(\theta)|^2\le \frac{2n^2 L^2}{K} \psi'(\theta)|\hat \mu^{2n}|(\theta)+ \frac{2n^4L^2}{K^3} \psi^2(\theta)
\psi'(\theta)|\hat \mu^{n-1}|(\theta)\, ,
\end{equation*}
Hence, arguing as in the proof of Proposition \ref{rittprop}, 
we see that $(\sigma_n)_{n\in \N}$ satisfies item $(ii)$ in Proposition 
\ref{ritt-ext}. 

%The fact that $(\sigma_n)_{n\in \N}$ satisfies items 
%$(ii)$ and $(iii)$ in Proposition \ref{genprop2} may be proved similarly. 

\medskip

It remains to prove the second part of the Proposition. 

\smallskip

Let $n\ge 1$. We have 
\begin{gather*}
(\delta_0-\sigma*\mu)*(\sigma*\mu)^{*n}= \sigma^{*n}*
\big[(\delta_0-\mu)* \mu^{*n}\big]+ \big[(\delta_0-\sigma)* 
\mu^{*(n+1)}\big]* \sigma^{*n}\, ,
\end{gather*}
which proves that $\sigma*\mu$ is Ritt. 

\medskip

Let $\alpha\in (0,1]$ and $n\ge 1$, and $\tau:=\alpha\mu+(1-\alpha 
)\sigma$. We have
\begin{gather*}
(\delta_0 -\tau)*\tau^{*n} =\\
\sum_{k=0}^n \dbinom{n}{k}\frac{\alpha^k (1-\alpha)^{n-k}}{k+1} 
\Big[\alpha (k+1)(\delta_0 -\mu)*\mu^{*k}+ (1-\alpha) (k+1)(\delta_0 
-\sigma )*\mu^{*k}\Big]*\sigma^{*(n-k)}\, .
\end{gather*}
Hence, 
$$
(n+1)\|(\delta_0 -\tau)*\tau^{*n}\|_{\ell^1} \le \frac{C}{\alpha} 
\sum_{k=0}^n \dbinom{n+1}{k+1} \alpha ^{k+1}(1-\alpha )^{(n+1)-(k+1)}
\le C\, ,
$$
and we see that $\tau$ is Ritt. \hfill $\square$

\smallskip

\subsection{Examples}

To exhibit examples we will make use of Proposition \ref{widder}. 
Hence we shall first exhibit completely monotone functions. 

\begin{lem}[Miller-Samko \cite{MS}]
Let $f,g\, :\, (0,+\infty)\to (0,+\infty)$ be infinitely differentiable functions functions such that $g'$ is completely monotone. 
\begin{itemize}
\item [$(i)$] If $f$ is completely monotone then $f\circ g$ is completely 
monotone either;
\item [$(ii)$] If $f'$ is completely monotone then $(f\circ g)'$ is completely monotone either.
\end{itemize}
%In particular if $f_1,\ldots , f_n\, :\, [0,+\infty)\to [0,+\infty)$ 
%admits completely monotone derivatives,  
\end{lem}
\noindent {\bf Proof.} 
Item $(i)$ is just Theorem 2 of \cite{MS}. Let us prove item $(ii)$. 
We have $(f\circ g)'=f'\circ g \,\times\,  g'$. By $(i)$, 
$f'\circ g$ is completely monotone. Then, $(f\circ g)'$ is completely 
monotone by Theorem 1 of \cite{MS}. \hfill $\square$

Define by induction $L_1(x)=L(x):=\log(1+x)$ and $L_{k+1}(x)=L (L_k(x))$ 
for every $x> 0$.

\begin{cor}
For every integer $k\ge 1$ and every real numbers $\alpha_1,\ldots ,\alpha_k \in [0,+\infty)$ and  $\alpha\in [0,+\infty)$ the function given by 
$$
f_{\alpha,\alpha_1,\ldots ,\alpha_k}(x)= \frac1{x^\alpha L_1(x)^{\alpha_1} 
\ldots L_k^{\alpha_k}(x)} \qquad \forall x\ge 0\, ,
$$
is completely monotone.
\end{cor}
\noindent {\bf Proof.} Obviously, $x\mapsto x^{-\alpha}$ is completely monotone. 
By $(ii)$ of the previous lemma $L_k$ admits a completely monotone derivative 
and then $L_k^{-\alpha_k}$ is completely monotone by $(i)$. The fact that 
$f_{\alpha,\alpha_1,\ldots ,\alpha_k}$ is also completely monotone then follows from Theorem 1 of \cite{MS}. \hfill $\square$

\medskip

\noindent {\bf Example 1.} Let $\mu$ be a probability measure supported on $\N$ 
such that $\mu(n)=cf_{\alpha,\alpha_1,\ldots ,\alpha_k}(n+1)$ for every $n\in 
\N$, where $\alpha_1,\ldots , \alpha_k\in [0,+\infty)$, $\alpha\in (1,2)$ and 
$c$ is a normalizing constant ensuring that we have a probability. Then, 
$\mu \in \HH \cap \RR$ 
Of course one may take $\alpha=1$ and $\alpha_1>1$, and so on... But 
for $\alpha=2$, $\mu$ does not even have BAR.
\medskip

It is more difficult to produce examples allowing negative $\alpha_k's$. 
One way to handle the difficulty is to proceed as in the proof of Proposition 
5.11 of \cite{CCL}.
\medskip

\noindent {\bf Example 2.} Our next example is a basic 
example of Ritt probability measures already considered by Dungey 
\cite{Dungey11} and Gomilko and Tomilov \cite{GT1}. Let 
$\gamma\in (0,1)$. We have a power series expansion 
$1-(1-t)^\gamma=\sum_{n\ge 1}a_n(\gamma)t^n$, $0\le t\le 1$. Notice that 
$\sum_{n\ge 1}a_n(\gamma)=1$ and $a_n(\gamma)\ge 0$ for every $n\ge 1$.  Define 
two probability measures $\tau$ and $\mu$ by setting for every $n\in \N$, $\mu(n)=a_{n+1}(\gamma)=\tau(n+1)$. so that 
$\tau=\delta_1* \mu$. Then, see for instance example 3.10a of \cite{GT1}, 
$\tau$ is a CM probability measure which has BAR. In particular, 
$\tau \in \HH \cap \RR$ and $\mu \in \HH\cap \RR$.

%\smallskip

%We emphasize here that the measure $\mu_2$ do not satisfy \eqref{cns1}.

 \section{Probability measures with a first moment}
 
 When $\mu$ has a first moment, a necessary condition 
 for the BAR property is that $\mu$ be centered, i.e. 
 $\sum_{n\in \Z} n \mu(n)=0$, see Proposition 1.9 of \cite{BJR94}. 
 
 \smallskip
 
 Hence we cannot consider probability measures $\mu$ supported by 
 $\N$ anymore. We shall consider the following situation. 
 
 \smallskip
 
 \begin{defn}
 We say that a probability measure $\mu$ on $\Z$ is CCM if it is supported on $\{-1\}\cup\N$ and if there exists 
  a finite positive measure $\nu$ on $[0,1]$, such that 
 $$
\int_0^1\frac{\nu(dt)}{(1-t)^2}=1\, . 
 $$  
 and 
 \begin{gather*}
 \mu(n):= \int_0^1t^n \nu(dt) \qquad \forall n\in \N;\\
 \mu(-1)= 1-\int_0^1\frac{\nu(dt)}{1-t}= 
 \int_0^1 \frac{t\nu(dt)}{(1-t)^2}\, .
 \end{gather*}
 
 \end{defn}
 \smallskip
 
 It is not hard to see that $\mu$ is indeed a probability measure and that it is centered.

 \subsection{Characterization of the BAR property}
 
% We shall not give the full details of  the proofs here, since many arguments are similar to the previous section. 

 \smallskip
 
 Let $\mu$ be a CCM probability measure on $\Z$ with representative measure $\nu$.
 
 \smallskip
 For every $\theta\in [-\pi,\pi]$, we have 
 \begin{equation*}
 \hat \mu(\theta)= \int_0^1\frac{1-2t +2t^2{\rm e}^{-i\theta}
 -t^2{\rm e}^{-2i\theta} }{(1-t)^2((1-t)^2+2t(1-\cos \theta))}
 \, \nu(dt)\, .
 \end{equation*}

In particular, 

\begin{equation}\label{real}
1-{\rm Re}\, \hat \mu(\theta)= (1-\cos \theta) \int_0^1\frac{2t(1-t\cos \theta)}{(1-t)^2((1-t)^2+2t(1-\cos \theta))}
 \, \nu(dt)\, ,
\end{equation}
and
\begin{equation}\label{imaginary}
{\rm Im}\, \hat \mu(\theta)= 2\sin \theta\, (1-\cos \theta) \int_0^1\frac{t^2}{(1-t)^2((1-t)^2+2t(1-\cos \theta))}
 \, \nu(dt)\,
\end{equation}

Consider the following condition on $\nu$: there exists $L>0$, such that 
for every $x\in [0,1)$,
\begin{equation}\label{condBARbis}
\frac{1}{1-x}\int_x^1\frac{t\nu(dt)}{(1-t)^2}\le L\int_0^x 
\frac{t\nu(dt)}{(1-t)^3}\, .
\end{equation}

Notice that if $\int_0^1\frac{\nu(dt)}{(1-t)^3}
<\infty$ (i.e. $\mu$ has a moment of order 2), condition \eqref{condBARbis} 
is automatically satisfied.  

\begin{prop}\label{propbis}
Let $\mu$ be a CCM probability measure on $\Z$ with representative measure $\nu$. Then, $\mu$ has BAR if and only if there exists $L>0$ 
such that $\nu$ satisfies \eqref{condBARbis}.
\end{prop}
\noindent {\bf Proof.} Assume \eqref{condBARbis}. Let us prove that 
$\mu$ satisfies \eqref{sector}. As noticed previously, it is enough to consider $\theta\in [-1/2,1/2]$. We have 
\begin{gather*}
|{\rm Im} \, \hat\mu(\theta)|\le C\Big(|\theta|^3  \int_0^{1-|\theta|} 
\frac{t}{(1-t)^4}\nu(dt) + |\theta|\int_{1-|\theta|}^1 
 \frac{t\nu(dt)}{(1-t)^2}\Big)\,. 
\end{gather*}
Using that $1-t\cos \theta\ge 1-t$, we see that 
$$
1-{\rm Re} \, \hat \mu(\theta)\ge \tilde C \theta^2 \int_0^{1-|\theta|}
 \frac{t\nu(dt)}{(1-t)^3}\, ,
$$
and \eqref{sector} holds, by \eqref{condBARbis}. 

\smallskip

Let us prove that if $\mu$ has BAR, then \eqref{condBARbis} holds. 
There exists $C>0$ such that for every $\theta\in [-1/2,1/2]$, 
\begin{equation}\label{bra}
|{\rm Im} \, \hat\mu(\theta)| \le C (1-{\rm Re} \, \hat \mu(\theta))\, .
\end{equation}
Let $\theta\in [-1/2,1/2]$ and $\alpha\in (0,1]$. We have
\begin{gather*}
|{\rm Im} \, \hat\mu(\theta)| \ge \frac{|\theta| }{4(1+\alpha^2)}
\int_{1-\alpha |\theta|}^1 
\frac{t\nu(dt)}{(1-t)^2}\, .
\end{gather*}
It is not hard to prove that there exists $C_\alpha,D>0$ such that 
\begin{gather*}
1-t\cos \theta \le C_\alpha (1-t) \qquad \forall t\in [0,1-\alpha|\theta|] \, ;\\
 1-t\cos \theta \le D \alpha |\theta| \qquad \forall t\in (1-\alpha|\theta|,1]\,. 
\end{gather*}
Hence, using \eqref{bra}, we infer that
$$
\frac{|\theta| }{4(1+\alpha^2)} \int_{1-\alpha |\theta|}^1 
\frac{t\nu(dt)}{(1-t)^2}\le C\Big( C_\alpha |\theta|^2\int_0^{1-|\theta|} 
\frac{t\nu(dt)}{(1-t)^3} + \alpha |\theta| \int_{1-\alpha |\theta|}^1 
\frac{t\nu(dt)}{(1-t)^2}\Big)\, .
$$
Taking $\alpha=1/(8C)$ gives the desired result. \hfill $\square$

\medskip

As before, we shall now characterize the BAR property in terms of 
the coefficients of $\mu$. 

\begin{prop}\label{propcaracbis}
Let $\mu$ be a CCM probability measure on $\Z$ with representative measure $\nu$. Then, $\nu$ satisfies \eqref{condBARbis} if and only if 
there exists $L>0$ such that 
\begin{equation}\label{condBARbis2}
n \sum_{k\ge n} k\mu(k) \le L \sum_{k=1}^n k^2\mu(k) \qquad \forall n\in \N\, .
\end{equation}
\end{prop}
\noindent {\bf Proof.} Assume \eqref{condBARbis}. Let $n\ge 2$, we have
\begin{gather*}
n \sum_{k\ge n} k\mu(k) \le \int_0^{1-1/n} \sum_{k\ge n} k^2 t^k \nu(dt)
+n \int_{1-1/n}^1 \frac{t\nu(dt)}{(1-t)^2}\\
\le (1+L) \int_{0}^{1-1/n}  \frac{t\nu(dt)}{(1-t)^3}\, .
\end{gather*}
Now, for every $1\le \ell \le  n-1$ and every $t\in [1-1/\ell,1-1/(\ell +1)]$, 
we have $\sum_{k=1}^n k^2 t^k\ge t \sum_{k=1}^\ell k^2 {\rm e}^{-1}\ge 
Ct/(1-t)^3$, where we used that $(1-1/m)^{m-1}$ decreases to ${\rm e}^{-1}$. 
Hence, \eqref{condBARbis2} holds.

\smallskip

Assume that \eqref{condBARbis2} holds. Let $\gamma\in (0,1]$ and $n\ge 2$. 
For every $t\in [1-1/n,1]$, since $\gamma\le 1$, we have 
\begin{gather*}
\sum_{k\ge \gamma n} kt^k  =\sum_{k\ge 0} (k+ n)t^{k+ n}
\ge \frac{t(1-1/n)^{n} }{(1-t)^2}\ge \frac{t({2\rm e})^{-1}}{(1-t)^2} \, .
\end{gather*}
Hence 
\begin{gather*}
n\int_{1-1/n}^1 \frac{t\nu(dt)}{(1-t)^2} \le 2{\rm e} n \sum_{k\ge \gamma n}k\mu(k)
\le \frac{2nL{\rm e}}{[\gamma n]} \sum_{ k=1}^{[ \gamma n]} k^2\mu(k)\\
\le \frac{2nL{\rm e}}{[\gamma n]} \Big(\int_0^{1-1/n} \frac{t\nu(dt)}{(1-t)^3} 
+ [\gamma n]^2 \int_{1-1/n} ^1 \frac{t\nu(dt)}{1-t}\Big)\, ,
\end{gather*}
and we conclude by taking $\gamma$ small enough \hfill $\square$

%\subsection{Ritt property in $\ell^1$}

\begin{theo}\label{theocompbis}
Let $\mu$ be a CCM probability measure on $\Z$ with representative measure $\nu$. Assume that 
$(\mu(n))_{n\in \N}$ satisfies  \eqref{condBARbis2}. Then, $\mu$ is Ritt and for every $m\in \N$,
there exists $C_m>0$ such that for every $f\in \ell^1(\Z)$, 
\begin{equation*}
\sup_{\lambda >0} \lambda \#\{k\in \Z\,:\, \sup_{n\ge 1} n^m|\mu^{*n}*(\delta_0-\mu)^{*m}*f(k)|\ge \lambda\}
\le C_m\|f\|_{\ell^1}\, .
\end{equation*}
\end{theo}
\noindent {\bf Proof.} It suffices to check that $\mu$ satisfies hypothesis \H\, 
and to apply Propositions \ref{genprop} and \ref{genprop2}. 

\medskip

To check the conditions we must estimate $\hat \mu$ and its derivatives.

Define
\begin{equation}\label{psy}
\psi(\theta)=\theta^2 \int_0^1 \frac{t\nu(dt)}{(1-t)((1-t)^2+\theta^2)}
=\frac{t\nu(dt)}{1-t}- \int_0^1\frac{t(1-t)\nu(dt)}
{(1-t)^2+\theta^2}\, .
\end{equation}
Then, 
\begin{equation}\label{psy'}
\psi'(\theta)=2\theta\int_0^1\frac{t(1-t)\nu(dt)}
{((1-t)^2+\theta^2)^2}\, .
\end{equation}
Hence for every $\theta\in [0,1/2]$, we have 
\begin{gather}
\label{psibound} \frac{\theta^2}{2}\int_0^{1-\theta}\frac{t\nu(dt)}{(1-t)^3}
+\frac12\int_{1-\theta}^1\frac{t\nu(dt)}{1-t}\le \psi(\theta)\le \theta^2\int_0^{1-\theta}\frac{t\nu(dt)}{(1-t)^3}
+\int_{1-\theta}^1\frac{t\nu(dt)}{1-t}\\
\nonumber {\theta}\int_0^{1-\theta}\frac{t\nu(dt)}{(1-t)^3}
+\frac1{\theta^3}\int_{1-\theta}^1t(1-t)\nu(dt)\le \psi'(\theta)\le 2{\theta}\int_0^{1-\theta}\frac{t\nu(dt)}{(1-t)^3}
+\frac2{\theta^3}\int_{1-\theta}^1t(1-t)\nu(dt)\, .
\end{gather}

In particular, using \eqref{condBARbis}, we see that \eqref{psi-cond} 
holds. 

Let us compute the derivatives of $\hat \mu$. We shall not give the full details here. Using \eqref{real}, we infer that

\begin{equation*}
1-{\rm Re}\, \hat\mu(\theta)=\int_0^1 \frac{1-t\cos \theta}{(1-t)^2}
\nu(dt) - \int_0^1\frac{1-t\cos \theta}{(1-t)^2+2t(1-\cos\theta)}
\nu(dt)\, ;
\end{equation*}

\begin{gather}\label{real'}
{\rm Re}\, \hat\mu'(\theta)= 
2\sin\theta (1-\cos\theta)\int_0^1 \frac{t^2 }{(1-t)^2((1-t)^2+2t(1-
\cos \theta))} \nu(dt)\\ 
\nonumber \qquad \qquad \qquad +\sin\theta \int_0^1\frac{2t(1-t\cos \theta)}
{((1-t)^2+2t(1-\cos\theta))^2}
\nu(dt)
 \\\nonumber =-\sin \theta \int_0^1 \frac{t}{(1-t)^2} 
\nu(dt) -\sin \theta \int_0^1 \frac{t(1-t^2)}{((1-t)^2+2t(1-
\cos \theta))^2} \nu(dt) \, ;
\end{gather}

and 

\begin{gather}\label{real''}
{\rm Re}\, \hat\mu''(\theta)= -\cos \theta \int_0^1 \frac{t}{(1-t)^2} 
\nu(dt) - \cos\theta \int_0^1 \frac{t(1-t^2)}{((1-t)^2+2t(1-
\cos \theta))^2} \nu(dt) \\
\nonumber \qquad \qquad \qquad +4\sin^2\theta \int_0^1 \frac{t^2
(1-t^2)}{((1-t)^2+2t(1-\cos \theta))^3} \nu(dt)\, .
\end{gather}

Using \eqref{imaginary}, we infer that

\begin{equation*}
{\rm Im }\, \hat \mu(\theta)=\sin \theta \int_0^1 \frac{t}{(1-t)^2} 
\nu(dt)-\sin\theta \int_0^1\frac{t}{(1-t)^2+2t(1-\cos\theta)}
\nu(dt)\, ;
\end{equation*} 

\begin{gather}
\label{imaginary'}{\rm Im }\, \hat \mu'(\theta)= 2\cos \theta (1-\cos \theta)
\int_0^1 \frac{t^2 }{(1-t)^2((1-t)^2+2t(1-
\cos \theta))} \nu(dt)
\\ \nonumber = \cos\theta \int_0^1 \frac{t}{(1-t)^2} 
\nu(dt)+ \int_0^1 \frac{2t^2-t(1+t^2)\cos \theta }{((1-t)^2+2t(1-
\cos \theta))^2} \nu(dt)\, ;
\end{gather}

and 

\begin{gather}
\label{imaginary''}{\rm Im }\, \hat \mu''(\theta)
=-\sin \theta \int_0^1 \frac{t}{(1-t)^2} 
\nu(dt)+\sin \theta \int_0^1 \frac{t(1+t^2)}{((1-t)^2+2t(1-
\cos \theta))^2} \nu(dt)\\
\nonumber \qquad \qquad \qquad -4\sin\theta \int_0^1 \frac{t(2t^2-t(1+t^2)\cos \theta )}{((1-t)^2+2t(1-
\cos \theta))^3} \nu(dt)
\end{gather}

We now derive the necessary estimates on $\hat \mu$ and its derivatives. 

Using \eqref{real}, we infer that 

\begin{claim}\label{claim1bis}
There exists $C>0$ such that  
$1-{\rm Re}\, \hat \mu (\theta)\ge C \theta^2\int_0^{1-\theta}\frac{t}{(1-t)^3}
\nu(dt)$, for every $\theta\in (0,1/2]$.
\end{claim}

Using \eqref{real'}, we infer that

\begin{claim}\label{claim2bis}
There exists $C>0$ such that  
$|{\rm Re}\, \hat \mu '(\theta)|\le C \theta\int_0^{1-\theta}\frac{t}{(1-t)^3}
\nu(dt) +\frac{C}{\theta^2}\int_{1-\theta}^1 t\nu(dt)$, for every $\theta\in (0,1/2]$.
\end{claim}

Using \eqref{imaginary'}, we infer that

\begin{claim}\label{claim3bis}
There exists $C>0$ such that  
$|{\rm Im}\, \hat \mu '(\theta)|\le C \theta^2\int_0^{1-\theta}\frac{t}{(1-t)^4}
\nu(dt) +C\int_{1-\theta}^1 \frac{t}{(1-t)^2}\nu(dt)$, for every $\theta\in (0,1/2]$.
\end{claim}

Using \eqref{real''}, we infer that

\begin{claim}\label{claim4bis}
There exists $C>0$ such that  
$|{\rm Re}\, \hat \mu ''(\theta)|\le C \int_0^{1-\theta}\frac{t}{(1-t)^3}
\nu(dt) +\frac{C}{\theta^4}\int_{1-\theta}^1 t(1-t)\nu(dt$, for every $\theta\in (0,1/2]$.
\end{claim}

Notice that there exists $\alpha>0$ such that for every $t\in [0,1]$ and every $\theta\in (0,1/2]$, 
$$|2t^2-t(1+t^2)\cos\theta)|=t|(1+t^2)(1-\cos \theta)-(1-t)^2|\le 
\alpha \max(\theta^2,(1-t)^2)\, .$$ 

Combining this estimate with \eqref{imaginary''}, we infer that

\begin{claim}\label{claim5bis}
There exists $C>0$ such that  
$|{\rm Im}\, \hat \mu ''(\theta)|\le C \theta\int_0^{1-\theta}\frac{t}{(1-t)^4}
\nu(dt) +\frac{C}{\theta^3}\int_{1-\theta}^1 t\nu(dt)$, for every $\theta\in (0,1/2]$.
\end{claim}

\medskip

We already saw that \eqref{psi-cond} holds. Let us prove that items $(i)-(iv)$ 
of Proposition \ref{genprop} hold. 

Item $(i)$ follows from Claim \ref{claim1bis} and \eqref{condBARbis} 
(see the proof of Theorem \ref{theocomp}. 

Item $(ii)$ follows from Claims \ref{claim2bis} and \ref{claim3bis} 
combined with \eqref{condBARbis} and \eqref{psibound}. 

Item $(iii)$ follows from item $(ii)$ combined with \eqref{psi-cond}. 

Item $(iv)$ follows from Claims \ref{claim4bis} and \ref{claim5bis} 
combined with \eqref{condBARbis} and \eqref{psibound}. \hfill $\square$

\medskip

\begin{prop}
Let $\tau$ be a  centered probability measure on $\Z$ such that $\sum_{n\in \Z}
|n|\tau(n) <\infty$. Assume moreover that 
there exists a CCM probability measure $\mu$ satisfying 
\eqref{condBARbis2} and that there exists  $a>0$ such that $\sum_{n\in \Z}n^2|\tau(n)-a\mu(n)|<\infty$. Then the conclusion of Theorem \ref{theocompbis} holds for $\tau$. 
\end{prop}
\noindent {\bf Proof.} We shall assume that $\sum_{n\in \Z} n^2 
\tau(n)=\infty$, otherwise, the result holds by Theorem \ref{centered}. In particular we must have $\sum_{n\in \Z}n^2 \mu(n)=\infty$ and by \eqref{psy} and \eqref{psy'} 

\begin{equation}\label{if}
\liminf_{\theta\in (0,\pi]}\psi(\theta)/\theta^2 =+\infty\quad \mbox{and}\quad  
\liminf_{\theta\in (0,\pi]}\psi'(\theta)/\theta =+\infty \qquad (\theta\to 0)\, . 
\end{equation}

It follows from the proof of Theorem \ref{theocompbis} 
that there exists an even function $\psi$ continuous on $[-\pi,\pi]$ and 
continuously differentiable on $(0,\pi]$, with $\psi(0)=0$ such that 
$\mu$ and $\psi$ satisfy the item $(i)-(iv)$ of hypothesis \H, 
for some $C,c>0$.

Since $\hat\tau=(\hat \tau -a \hat \nu )+a\hat \mu$ is clearly twice differentiable on $(0,\pi]$, the proposition will be proved if we can show that the items $(i)-(iv)$ of hypothesis \H \,
hold with $\tau$ in place of $\mu$ with the same $\psi$, but for possibly different $C,c>0$. 

\smallskip

We already saw that $\tau$ must be strictly aperiodic. Hence $|\hat \mu|<1$ 
on $(0,\pi]$. In particular, to prove item $(i)$ it suffices to consider 
$\theta\in (0,\eta]$ for some small enough $\eta>0$. 

For every $\theta\in (0,\pi]$, we have 
$$
\hat\tau (\theta)= \sum_{n\in \Z} (\tau(n)-a\mu(n))({\rm e}^{in\theta} 
-1) + \big[1-a+a\sum_{n\in\Z}\mu(n) {\rm e}^{in\theta}\big]:= 
\chi(\theta)+ \phi(\theta)\, .
$$
Using that $\sum_{n\in \Z}n^2|\tau(n)-a\mu(n)|<\infty$ and that 
$\sum_{n\in \Z}n(\tau(n)-a\mu(n))=0$, we see that $\lim_{\theta \to 0,
\theta\neq 0} \chi(\theta)/\theta^2=\chi''(0)$ exists. In particular, 
$\lim_{\theta \to 0,
\theta\neq 0} \chi(\theta)/\psi(\theta)=0$. 

Now, we have 
$$
|\phi(\theta)|^2=(1-a+a{\rm Re}\, \hat \mu(\theta))^2 +a^2({\rm 
Im}\, \hat \mu(\theta) )^2\, =1-2a(1-{\rm Re}\, \hat\mu(\theta)) 
+a^2(1-{\rm Re}\, \hat\mu(\theta)) ^2+a^2 ({\rm 
Im}\, \hat \mu(\theta) )^2.
$$
Hence, using Claim \ref{claim1bis}, we infer that there exists $\eta >0$ such that for every 
$\theta\in (0,\eta]$, $|\hat \tau (\theta)|\le 1-\delta\psi(\theta)$, 
for some $\delta>0$.

\smallskip

The proofs of item $(ii)-(iii)$ are similar (but simpler) hence we leave them to 
the reader. \hfill $\square$

\medskip

\noindent {\bf Example 3.} Let $\alpha\in (2,+\infty)$ and $\alpha_1,\ldots 
,\alpha_k\ge 0$. Let $\mu$ be a probability on $\Z$ such that 
$\sum_{n\in \Z} |n|\mu(n)<\infty$, $\sum_{n\in \Z} n\mu(n)=0$ and 
$\sum_{n\in \Z}n^2|\mu(n)-a f_{\alpha,\alpha_1,\ldots , \alpha_k}(n+1)|<\infty$, 
for some $a>0$, where we extended $f_{\alpha,\alpha_1,\ldots , \alpha_k}$ 
to $\Z^-$ by setting, $f_{\alpha,\alpha_1,\ldots , \alpha_k}(-n)=0$ 
for every $n\in \N$.

\section{Symmetric probability measures} \label{symmetric}

In this section we consider symmetric probability measures. If $\mu$ 
is symmetric (i.e. $\check \mu=\mu)$, then  
$\hat \mu$ is real valued, hence has BAR. It is known that 
if moreover $(\mu(n))_{n\in \N}$ is  non-increasing then 
\eqref{weakinemax0} holds. However we are not aware of any result 
concerning the Ritt property or \eqref{weakinemax} with 
$m\ge 1$. 

We shall again investigate the situation where we have  completely 
monotone coefficients.  To be more precise we consider the following situation.

\begin{defn}
We say that a probability measure $\mu$ on $\Z$ is SCM if it is symmetric and 
if there exists a finite positive measure on $[0,1]$ such that 
\begin{gather*}
\int_0^1 \frac{1}{1-t}\nu(dt)=1/2 \quad , \quad \mu(0)=2
\int_0^1 \nu(dt);\\
\mu(n)= \int_0^1t^n\nu(dt) \qquad \forall n\ge 1\, .
\end{gather*}
\end{defn}

Let $\mu$ be an SCM probability measure on $\Z$ with repesentative measure $\nu$. 
Define another measure on $\Z$, supported on $\N$, by setting 
\begin{gather*}
 \quad \mu_1(0)=2
\int_0^1 \nu(dt);\\
\mu_1(n)= 2\int_0^1t^n\nu(dt) \qquad \forall n\ge 1\, .
\end{gather*}

Then $\mu_1$ is a probability measure, $(\mu_1(n))_{n\in \N}$ is completely monotone and  $\mu=\frac12(\check \mu_1+\mu_1)$. In particular, it follows from 
Proposition \ref{ext} and Theorem \ref{theocomp}, that 
$\mu$ satisfies hypothesis \H\,  as soon as $\mu$ satisfies \eqref{condBAR2}. 
The fact that $\mu$ is Ritt when it satisfies \eqref{condBAR2} may be proved similarly (but more easily). 

\medskip

We could use a similar argument based on Theorem \ref{theocompbis}. 
However, doing so, we would miss some symmetric probability measures 
satisfying hypothesis \H.

Let us explain how to be more precise. Let $\mu$ be a SCM probability measure. 

It follows from previous computations that, for every $\theta\in \R$,

$$
1-\hat \mu(\theta)= 1- {\rm Re\, }\hat \mu(\theta) = 
\int_0^1 \frac{t(1-\cos\theta)}{(1-t)((1-t)^2+2t(1-\cos \theta))}\,\nu(dt)\,  .
$$

Consider the following condition on $\nu$: there exists $L>0$ such that 
for every $x\in [0,1)$,
\begin{equation}\label{condBARter}
\int_x^1 \frac{t}{1-t} \nu(dt) \le L (1-x)^2 \int_0^x \frac{t}{(1-t)^3}\nu(dt)
\, .
\end{equation}
This condition can be proved to be equivalent to the following one: there exists 
$D>0$ such that for every $n\ge 1$, 
\begin{equation}\label{condBARter2}
n^2 \sum_{k\ge n}\mu(k) \le L \sum_{k=1}^n k^2\mu(k)\, .
\end{equation}

One can prove that if \eqref{condBARter} holds, then $\mu$ satisfies hypothesis 
\Ht\, with $\psi$ given by 
\begin{equation*}
\psi(\theta)= \theta^2 \int_0^1\frac{t}{(1-t)(1-t+|\theta|)^2} 
\nu(dt) \qquad \forall \theta \in [-\pi,\pi]-\{0\}\, .
\end{equation*} 
Notice that $\psi'(\theta)= 2\theta \int_0^1\frac{t}{(1-t+\theta)^3} \nu(dt)$, 
for every $\theta\in (0,\pi]$. 

\medskip

Then, one can prove that a SCM probability measure satisfying 
\eqref{condBARter2} is Ritt and satisfies \eqref{weakinemax} for every $m\in \N$ and some $C_m>0$. 

In particular, we have the following.

\begin{theo}\label{theosymmetric}
Let $\mu$ be a SCM probability measure such that $(\mu(n))_{n\in \N}$ satisfies 
 either \eqref{condBAR2} or \eqref{condBARter2}. Then, $\mu$ is Ritt 
 and satisfies \eqref{weakinemax} for every $m\in \N$ and some $C_m>0$.
\end{theo}

\medskip

\noindent {\bf Example 4.} Let $\alpha>1$ and $\alpha_1, \ldots, \alpha_k\ge 0$. 
Let $\mu$ be a symmetric probability measure defined by 
$\mu(0)=2cf_{\alpha,\alpha_1,\ldots ,\alpha_k}(1)$ and for every $n\ge 1$ 
$\mu(n)=cf_{\alpha,\alpha_1,\ldots ,\alpha_k}(n+1)$, where $c$ is a normalizing sequence ensuring that $\mu$ is a probability. Then, $\mu$ is a SCM probability measure for which the above theorem apply.

\section{Discussion and open questions} 

- Most of the examples of (strictly aperiodic) probability measures on $\Z$ that have BAR are known to be Ritt. We do not believe that the BAR property and the Ritt property are equivalent, but one has to find a counterexample. This problem 
was also formulated by Dungey  \cite{Dungey11} (see his remarks page 1729). 

\medskip

- One may wonder whether, in the symmetric case, the condition "$(\mu(n))_{n
\in \N}$ is non-increasing" is sufficient for the Ritt property or for weak type maximal inequalities \eqref{bc}, since it is sufficient 
for the weak type maximal inequality \eqref{weakL1}. At least, for a SCM 
probability measure on $\Z$, can one remove the conditions \eqref{condBAR2} and \eqref{condBARter2} from Theorem \ref{theosymmetric} ? 

\medskip

- Let $\mu$ be a probability measure  on $\Z$. Let $f\in \ell^p(\Z)$, $p\ge 1$. 
Consider the square function defined by $s_\mu(f)(k):=\Big(\sum_{n\ge 1} n 
\big((\mu^{*n}-\mu^{*(n+1)})*f(k)\big)^2\Big)^{1/2}$. 
Assume that $\mu$ has BAR. When $p>1$, it follows from the work of Le Merdy 
and Xu \cite{LX12} that there exists $C_p>0$ such that for every 
$f\in \ell^p(\Z)$, $\|s_\mu(f)\|_p\le C_p \|f\|_p$, i.e. $s(f)$ satifies a strong 
$p-p$ inequality. A natural question is whether $s(f)$ satisfies 
a weak $1-1$ inequality.

\end{document}